\documentclass{article}

\usepackage{microtype}
\usepackage{graphicx}
\usepackage[caption=false,font=normalsize,labelfont=sf,textfont=sf]{subfig}
\usepackage{booktabs} % for professional tables

\usepackage{hyperref}

\usepackage[accepted]{icml2024}

\usepackage{amsmath}
\usepackage{amssymb}
\usepackage{mathtools}
\usepackage{amsthm}
\usepackage{mathrsfs}
\usepackage{multirow}
\usepackage{algorithm}
\usepackage{algorithmic}
\usepackage{bm}
\usepackage{url}
\usepackage{float}
\usepackage{enumitem}

\usepackage[capitalize,noabbrev]{cleveref}

\theoremstyle{plain}
\newtheorem{theorem}{Theorem}[section]
\newtheorem{proposition}[theorem]{Proposition}
\newtheorem{lemma}[theorem]{Lemma}

\theoremstyle{definition}
\newtheorem{definition}[theorem]{Definition}

\theoremstyle{remark}

\usepackage[textsize=tiny]{todonotes}

\icmltitlerunning{Autonomous Sparse Mean-CVaR Portfolio Optimization}

\begin{document}

\def \bbE {\mathbb E}
\def \bbR {\mathbb R}
\def \bbN {\mathbb N}
\def \bbZ {\mathbb Z}

\def \mA {\mathcal{A}}
\def \mD {\mathcal{D}}
\def \mF {\mathcal{F}}
\def \mG {\mathcal{G}}
\def \mI {\mathcal{I}}
\def \mJ {\mathcal{J}}
\def \mL {\mathcal{L}}
\def \mM {\mathcal{M}}
\def \mN {\mathcal{N}}
\def \mO {\mathcal{O}}
\def \mP {\mathcal{P}}
\def \mQ {\mathcal{Q}}
\def \mS {\mathcal{S}}
\def \mT {\mathcal{T}}
\def \mU {\mathcal{U}}
\def \mZ {\mathcal{Z}}

\def \ba {\bm{a}}
\def \bb {\bm{b}}
\def \bc {\bm{c}}
\def \bd {\bm{d}}
\def \bg {\bm{g}}
\def \bh {\bm{h}}
\def \bp {\bm{p}}
\def \bq {\bm{q}}
\def \bx {\bm{x}}
\def \by {\bm{y}}
\def \bz {\bm{z}}
\def \bw {\bm{w}}
\def \bu {\bm{u}}
\def \bv {\bm{v}}
\def \br {\bm{r}}
\def \bs {\bm{s}}
\def \bR {\bm{R}}
\def \bS {\bm{S}}
\def \bI {\bm{I}}
\def \bA {\bm{A}}
\def \bB {\bm{B}}
\def \bC {\bm{C}}
\def \bD {\bm{D}}
\def \bE {\bm{E}}
\def \bF {\bm{F}}
\def \bG {\bm{G}}
\def \bH {\bm{H}}
\def \bL {\bm{L}}
\def \bP {\bm{P}}
\def \bQ {\bm{Q}}
\def \bR {\bm{R}}
\def \bS {\bm{S}}
\def \bU {\bm{U}}
\def \bV {\bm{V}}
\def \bW {\bm{W}}
\def \bX {\bm{X}}
\def \bY {\bm{Y}}

\def \blambda {\bm{\lambda}}
\def \bLambda {\bm{\Lambda}}

\def \bbRn {\bbR^{n}}
\def \bbRm {\bbR^{m}}
\def \bbRN {\bbR^{N}}
\def \bbRNN {\bbR^{N\times N}}
\def \bbRM {\bbR^{M}}
\def \bbRT {\bbR^{T}}
\def \bbNM {\bbN_M}
\def \bbNk {\bbN_k}
\def \bbNn {\bbN_n}
\def \bbNm {\bbN_m}

\def \tbQ {\tilde{\bQ}}
\def \tbq {\widetilde{\bq}}
\def \tbI {\widetilde{\bI}}

\def \st {\text{s.t.}}
\def \prox {\mathrm{prox}}
\def \diag {\mathrm{diag}}
\def \VaR {\mathrm{VaR}}
\def \CVaR {\mathrm{CVaR}}

\def \Argmax {\mathop{\rm Arg\,max}}
\def \Argmin {\mathop{\rm Arg\,min}}

\def \argmax {\mathop{\rm arg\,max}}
\def \argmin {\mathop{\rm arg\,min}}

\newcommand\leqs{\leqslant}
\newcommand\geqs{\geqslant}

\twocolumn[
\icmltitle{Autonomous Sparse Mean-CVaR Portfolio Optimization}



\icmlsetsymbol{equal}{*}

\begin{icmlauthorlist}
\icmlauthor{Yizun Lin}{affil1}
\icmlauthor{Yangyu Zhang}{affil2}
\icmlauthor{Zhao-Rong Lai}{affil1}
\icmlauthor{Cheng Li}{affil1}
\end{icmlauthorlist}

\icmlaffiliation{affil1}{Department of Mathematics, College of Information Science and Technology, Jinan University, Guangzhou, China}
\icmlaffiliation{affil2}{Jinan University-University of Birmingham Joint Institute, Jinan University, Guangzhou, China}

\icmlcorrespondingauthor{Zhao-Rong Lai}{laizhr@jnu.edu.cn}

\icmlkeywords{Mean-CVaR, $\ell_0$ constraint, tailed approximation, proximal alternating linearized minimization}

\vskip 0.3in
]



\printAffiliationsAndNotice{This is a preprint version of an accepted paper for ICML 2024. Contents may be changed in its final version.}  % leave blank if no need to mention equal contribution

\begin{abstract}
The $\ell_0$-constrained mean-CVaR model poses a significant challenge due to its NP-hard nature, typically tackled through combinatorial methods characterized by high computational demands. From a markedly different perspective, we propose an innovative autonomous sparse mean-CVaR portfolio model, capable of approximating the original $\ell_0$-constrained mean-CVaR model with arbitrary accuracy. The core idea is to convert the $\ell_0$ constraint into an indicator function and subsequently handle it through a tailed approximation. We then propose a proximal alternating linearized minimization algorithm, coupled with a nested fixed-point proximity algorithm (both convergent), to iteratively solve the model. Autonomy in sparsity refers to retaining a significant portion of assets within the selected asset pool during adjustments in pool size. Consequently, our framework offers a theoretically guaranteed approximation of the $\ell_0$-constrained mean-CVaR model, improving computational efficiency while providing a robust asset selection scheme.
\end{abstract}

\section{Introduction}
Value-at-risk (VaR, \citealt{valueatrisk}) is a widely-used downside-risk metric in finance that assesses the extreme percentage loss in a disadvantageous situation. Since it fails to satisfy the subadditivity property \citep{varsubadd}, it is improved to its tail statistic named conditional value-at-risk (CVaR, \citealt{convalueatrisk}). CVaR not only has the good theoretical property of coherence \citep{varsubadd,meancvar}, but also accords with dual stochastic dominance \citep{ssdorigin}. Moreover, CVaR is more tractable than VaR from the perspective of optimization. It can serve as the risk metric in Markowitz's mean-risk criterion \citep{PS} which balances the expected return and the risk of a portfolio. The original mean-CVaR model with linear constraints can be efficiently solved by linear programming \citep{convalueatrisk}.

In practical portfolio optimization (PO), it is desirable to constrain the scale of selected assets while keeping an appropriate investing performance, in order to reduce transaction cost and managerial burden. This can be done via the managerial approach, such as the endowment model \citep{POmasci2}, market crashes \citep{POmasci3}, and the revenue-driven resource allocation \citep{POmasci1}. However, these methods still demand professional knowledge in management and finance. Recently, the arising development in machine learning methods for management science provides an alternative way for PO. For example, \cite{mlpo} propose the performance-based regularization with a data-driven strategy to improve the out-of-sample performance in PO. \cite{sparsepo} adopt the $\ell_1$ regularization to construct a sparse and stable Markowitz portfolio (SSMP). \cite{SSPO} and \cite{SSPOl0} propose short-term sparse portfolio optimization models with $\ell_1$ and $\ell_0$ regularizations, respectively. \cite{SPOLC} propose a sparse structure for the covariance estimation based on the spectral decomposition.

Among all the sparsity approaches, the $\ell_0$ constraint can exactly control the number of selected assets. However, the $\ell_0$ constraint generally leads to an NP-hard problem that has a very high computational complexity. A general approach to exactly solve the $\ell_0$-constrained problems is the combinatorial approach, such as the mixed-integer optimization \citep{MIO} and the branch-and-bound procedure \citep{cutplanecvar}. For example, \cite{cutplanePO} propose a cutting-plane approach with some heuristics to exactly solve an $\ell_0$-constrained mean-variance model, which saves some computational time according to the experimental results. \cite{cutplanecvar} further adopt this approach to exactly solve an $\ell_0$-constrained mean-CVaR model. It is essentially a branch-and-bound algorithm: for each visit of a node (a combinatorial case), it conducts one ordinary mean-CVaR optimization. Since its experimental results show that the number of visited nodes ranges from about 1,300 to more than 77,000, the computational complexity is still high in practical portfolio management. Moreover, \citep{cutplanecvar} requires the commercial optimizer Gurobi\footnote{\url{https://www.gurobi.com/}} to deploy this cutting-plane approach. We try to run the program, but it prompts an error that the models are too large in our experiments. In summary, the existing exactly-solving approaches have high computational complexity and rely on the external commercial optimizer, thus economic and practical approaches still need to be developed.
%\href{https://www.gurobi.com/}{Gurobi}

Instead of exactly solving the $\ell_0$-constrained mean-CVaR model, we propose an autonomous sparse mean-CVaR (ASMCVaR) PO model that can approximate it to arbitrary accuracy. The core idea is to convert the $\ell_0$ constraint into an indicator function and approximate it to arbitrary accuracy through a tailed approximation approach (Definition \ref{def:tailed}). By this way, we can replace the nonconvex and complex geometric structure of $\ell_0$ constraint by a simple and more tractable tailed approximation. To solve the ASMCVaR model, we develop a Proximal Alternating Linearized Minimization Algorithm (PALM) with a nested Fixed-Point Proximity Algorithm (FPPA) that can handle all the constraints and guarantee convergence. It saves much computational complexity and keeps a certain accuracy at the same time. ASMCVaR can keep a large proportion of assets in the selected asset pool when adjusting the pool size, which reduces managerial burden and delivers convenience. Our main contributions can be summarized as follows.
\begin{itemize}[topsep=0pt, itemsep=0pt, parsep=0pt]
\item[1.] We propose an ASMCVaR model and prove that it can approximate the $\ell_0$-constrained mean-CVaR model to arbitrary accuracy based on item 2.
\item[2.] We propose a tailed approximation for the indicator function of the $\ell_0$ constraint. We also find that this approximation can be arbitrarily accurate as the positive approximation parameter tends to $0$. See Definition \ref{def:tailed} and Figure \ref{fig:tail} for more detail.
\item[3.] We propose a convergent PALM algorithm with a nested FPPA to iteratively solve ASMCVaR. It does not require any external or commercial optimizers.
\item[4.] ASMCVaR can keep a large proportion of assets in the selected asset pool when adjusting the pool size, thus the selected assets are stable. It is beneficial to practical portfolio management when the manager needs to adjust the asset pool size.
\end{itemize}

\section{Preliminaries and Related Works}
\label{sec:relate}
We present some preliminaries on CVaR, and then introduce some approximate or exact solutions to the $\ell_0$-constrained mean-CVaR model.

\subsection{VaR and CVaR}
VaR measures the percentage loss of portfolio return in a bad scenario and with a certain confidence level $c$:
\begin{equation*}
\label{eq:var}
\mathrm{Pr}\{\VaR\leqs-r_{s}\} = 1 - c,
\end{equation*}
where $\mathrm{Pr}\{\cdot\}$ denotes the probability of an event and $r_{s}$ denotes the rate of return of a portfolio (representing some kind of investing strategy). Specifically, if $\VaR(c=95\%)=0.1$, then this investing strategy may suffer a loss greater than $10\%$ with a probability of $5\%$.

CVaR further averages the values greater than the VaR at a given confidence level $c$:
\begin{equation*}
\label{eq:cvar}
\CVaR:=\min\limits_{\tau\in\bbR}\left\{\tau+\frac{1}{1-c}\mathbb{E}\left(\max\{-r_{s}-\tau,0\}\right)\right\},
\end{equation*}
where $\mathbb{E}(\cdot)$ denotes the expectation of a random variable.

\subsection{Original Mean-CVaR Model}
In a real-world investing environment, the CVaR can be formulated as a linear program. Let $\bR:=[\br_{1}^{\top};\cdots;\br_{T}^{\top}] \in \bbR^{T \times N}$ be the sample return matrix of $N$ assets and $T$ observations gathered from a financial market. Accordingly, a portfolio (investing strategy) can be seen as an $N$-dimensional vector $\bw$. Then the sample returns of this portfolio are $\{\br_t^{\top}\bw\}_{t=1}^T$. Using the sample average estimator for CVaR yields
\begin{align}
\label{eq:cvaroptim}
&\min\limits_{\bw\in\bbRN,\ \tau\in\bbR,\ \bz\in\bbRT}\left\{\tau+\frac{1}{(1-c)T}\bm{1}^{\top}_{T}\bz\right\}\\
\label{eq:cvarineq}
&\st\quad \bz\geqs-\bR\bw-\tau{\bm1}_T,\quad \bz\geqs\bm{0}_{T},\\
\label{eq:cvarineq1}
&\ \bw\in \Delta:=\left\{\bw\in\bbR^N:  \bw\geqs \bm{0}_{N}, \bw^\top \bm{1}_{N}=1\right\},   \\
\label{eq:cvarineq2}
&\ (\bw^\top\hat{\bm{\mu}}=\rho),
\end{align}
where $\bz$ is the auxiliary variable to deal with the positive part in the expectation. $\bm{1}_{T}$ or $\bm{0}_{T}$ is a vector of $1$ or $0$ with dimensionality shown in the subscript, respectively. The inequality for vectors dominates each dimension. In \eqref{eq:cvarineq}, inequality constraints for $\bz$ can be used instead of equality constraints because $\bz$ is greater than $\bm{0}_{T}$ and the objective is a minimization. \eqref{eq:cvarineq1} is the simplex constraint that consists of the no-short-selling constraint and the self-financing constraint, respectively. Such constraints ensure feasibility of the portfolio in real-world investment. \eqref{eq:cvarineq2} is an optional constraint for the expected portfolio return, where $\hat{\bm{\mu}}=\frac{1}{T}\bR^\top \bm{1}_{T}$ is the sample asset return vector. \eqref{eq:cvaroptim}$\sim$\eqref{eq:cvarineq2} constitute the original mean-CVaR model \citep{convalueatrisk}. Since all the constraints are convex and linear, this model can be efficiently solved by the simplex method.

\subsection{Mean-CVaR Model with $\ell_1$ Regularization}
\label{sec:l1meancvar}
A widely adopted approach for sparsity is to impose $\ell_1$ regularization. This method is favored for its compatibility with widely recognized solving techniques, such as the linear programming, the Lasso \citep{LASSO} and the LARS \citep{LARS} algorithms. \cite{meancvarl1} propose a mean-CVaR model with a re-weighted $\ell_1$ regularization as follows:
\begin{equation}
\label{model:MeanCVaRl1}
\begin{split}
&\min\limits_{\bw\in\bbRN,\tau\in\bbR,\bz\in\bbRT}\left\{\tau+\frac{1}{(1-c)T}\bm{1}^{\top}_{T}\bz+\blambda^{\top}\bw\right\}\\
&\st\ \ \bz{\geqs-}\bR\bw{-}\tau{\bm1}_T,\bz\geqs\bm{0}_{T},\bw\in\Delta,(\bw^\top\hat{\bm{\mu}}=\rho),
\end{split}
\end{equation}
where $\blambda\in\{\bx\in\bbRN: \bx\geqs \bm{0}_{N}\}$ is the weight vector for the $\ell_1$ regularization. Since $\bw\geqs \bm{0}_{N}$, the re-weighted $\ell_1$ regularization in \eqref{model:MeanCVaRl1} can be simplified: $\sum_{i=1}^N\lambda_i |w_i|=\sum_{i=1}^N\lambda_i w_i=\blambda^{\top}\bw$. This formulation can also be efficiently solved by linear programming. In each outer-iteration, the algorithm updates $\blambda$ by some rules regarding $\bw$, for example, $\lambda_i^{k+1}=1/(|w_i^{k}|+\epsilon)$. Then in the inner-iteration, \eqref{model:MeanCVaRl1} is solved with fixed $\blambda^{k+1}$ to obtain $\bw^{k+1}$. This approach enjoys a low computational cost, but it cannot exactly control the number of selected assets. One can just empirically tune $\blambda$ to get a rough number of selected assets. Besides, it is also intractable to assess the accuracy of approximation of $\ell_1$ regularization to the $\ell_0$ constraint.

\subsection{$\ell_0$-constrained Mean-CVaR Model}
\label{sec:sparsemeancvar}
To impose exact sparsity on the portfolio, the $\ell_0$ constraint should be used:
\begin{align}
\label{model:MeanCVaR}
&\min\limits_{\bw\in\bbRN,\ \tau\in\bbR,\ \bz\in\bbRT}\left\{\tau+\frac{1}{(1-c)T}\bm{1}^{\top}_{T}\bz\right\}\\
%\label{model:MeanCVaR2}
\notag&\st\ \bz\geqs-\bR\bw-\tau{\bm1}_T,\ \bz\geqs\bm{0}_{T},\ \bw\in\Delta,\ (\bw^\top\hat{\bm{\mu}}=\rho),\\
\label{model:MeanCVaRl0}
&\hspace{1.5em}\|\bw\|_{0}\leqs m,
\end{align}
where $\Delta$ is defined in \eqref{eq:cvarineq1}, and $\|\bw\|_{0}$ represents the number of nonzero components in $\bw$. The $\ell_0$ constraint in \eqref{model:MeanCVaRl0} is also called the $m$-sparse constraint, since $m$ controls the number of selected assets. The above $\ell_0$-constrained mean-CVaR problem is NP-hard. To solve it by the combinatorial approach, \cite{cutplanecvar} introduce an auxiliary variable $\by\in \mZ_N^m:=\{\by\in \{0,1\}^N: \bm{1}_{N}^{\top}\by\leqs m\}$, then \eqref{model:MeanCVaRl0} is equivalent to the logical constraint
\begin{equation}
\setlength{\abovedisplayskip}{3pt}
\setlength{\belowdisplayskip}{3pt}
\label{model:MeanCVaRl0logic}
  \by\in \mZ_N^m, \quad y_i=0 \Rightarrow w_i=0.
\end{equation}
Let $\bY:=\diag(\by)$ be the diagonalized matrix for $\by$. Then model \eqref{model:MeanCVaR} can be reformulated as a bilevel optimization problem
\begin{align}
\setlength{\abovedisplayskip}{4pt}
\setlength{\belowdisplayskip}{4pt}
\label{model:uplevel}
&\text{(upper-level)}\qquad  \min_{\by\in \mZ_N^m}  f(\by),\\
\label{model:lowlevel1}
&\text{(lower-level)}\nonumber\\
&f(\by)=\min_{\bw\in\bbRN,\tau\in\bbR,\bz\in\bbRT} \ \left\{\frac{1}{2\gamma} \bw^\top\bw{+}\tau{+}\frac{1}{(1{-}c)T}\bm{1}^{\top}_{T}\bz\right\}\\
\label{model:lowlevel2}
 & \st\ \bz{\geqs-}\bR\bY\bw{-}\tau{\bm1}_T,\bz{\geqs}\bm{0}_{T},\bY\bw{\in}\Delta,(\bw^\top\bY\hat{\bm{\mu}}{=}\rho),
\end{align}
where the logical constraint \eqref{model:MeanCVaRl0logic} is embedded in $\bY\bw$. Note that $\bw=\bY\bw$ after minimization \eqref{model:lowlevel1}, thus $(\bY\bw)^\top\bY\bw$ is replaced by $\bw^\top\bw$ \citep{cutplanePO}. As $\gamma\rightarrow +\infty$, this bilevel optimization can approximate the $\ell_0$-constrained mean-CVaR optimization to arbitrary accuracy.

In the upper-level optimization, \cite{cutplanecvar} solve the following problem
\begin{equation}
\setlength{\abovedisplayskip}{4pt}
\setlength{\belowdisplayskip}{4pt}
\label{model:uplevel1}
  \min_{(\by,\theta)\in \mZ_N^m\times \bbR} \theta \quad \st \quad (\by,\theta)\in \mF^k\subseteq\mZ_N^m\times \bbR
\end{equation}
to obtain the node $(\by^k,\theta^k)$ of the $k$-th visit, where $\mF^k$ denotes the $k$-th feasible region for $(\by,\theta)$. $\theta$ indicates any possible value for $f(\by)$ in the current iteration, and its infimum  provides a lower bound for $f(\by)$. To initialize, one can just solve a continuously-relaxed and lower-bounded version of \eqref{model:MeanCVaR}$\sim$\eqref{model:MeanCVaRl0} to obtain $\theta^1$ for all $\by\in \mZ_N^m$ and pick up any $\by^1\in \mZ_N^m$. In the lower-level optimization, \eqref{model:lowlevel1}$\sim$\eqref{model:lowlevel2} is a convex quadratic optimization for which efficient algorithms exist. If it is infeasible, then $\mF^{k+1}$ is obtained by removing $(\by^k,\theta^k)$ from $\mF^k$; or else it can be solved to obtain $f(\by^k)$. Since $f(\by)$ is convex in $\by\in[0,1]^N$, the set
\begin{equation}
\setlength{\abovedisplayskip}{4pt}
\setlength{\belowdisplayskip}{4pt}
\label{model:Fkplus1}
\{(\by,\theta)\in \mZ_N^m\times \bbR: \theta\geqs f(\by^k)+\bg_{\by^k}^\top(\by-\by^k)\}
\end{equation}
provides new lower bounds for all the unvisited cases $\by\in \mZ_N^m$, where $\bg_{\by^k}$ is a subgradient of $f$ at $\by^k$. Intersecting \eqref{model:Fkplus1} and $\mF^{k}$ yields $\mF^{k+1}$. Either feasible or infeasible, the lower-level optimization narrows down the searching region from $\mF^{k}$ to $\mF^{k+1}$. Then in the next visit (iteration), \eqref{model:uplevel1} is solved in $\mF^{k+1}$. During this procedure, $LB^k:=\theta^k$ serves as the current uniform lower bound for all the unvisited cases, while $UB^k:=\min_{1\leqs j\leqs k} f(\by^j)$ serves as the current uniform upper bound for all the visited cases. If $(UB^k-LB^k<\epsilon)$ for a sufficiently small positive $\epsilon$, then $\{\by^j: 1\leqs j\leqs k, f(\by^j)=UB^k\}$ can be considered as the $\epsilon$-optimal solution set. In fact, this is a branch-and-bound strategy.

Though the above bilevel optimization takes the exact solving objective as the initial intention, it involves two approximations: one is the term $\frac{1}{2\gamma} \bw^\top\bw$ in \eqref{model:lowlevel1}, the other is the $\epsilon$ optimality of the gap $(UB^k-LB^k)$. Moreover, this branch-and-bound approach may visit tens of thousands of combinatorial cases to achieve satisfactory optimality. In addition, external commercial optimizers are also required to deploy the whole algorithm. The program is unable to run in our experiments because the models are too large. To develop economic and practical methods, we propose to approximate the $\ell_0$-constrained mean-CVaR model from a very different perspective based on a tailed approximation.

\section{Autonomous Sparse Mean-CVaR Portfolio Optimization}
In this section, we propose an Autonomous Sparse Mean-CVaR (ASMCVaR) model and develop a Proximal Alternating Linearized Minimization (PALM) algorithm to solve this model.

\subsection{Autonomous Sparse Mean-CVaR Model}
\label{sec:ASMCVaR}
Before introducing the ASMCVaR model, we revise the $\ell_0$-constrained Mean-CVaR model and reformulate it as a nonconvex three-term optimization model. Note that there is an equality constraint $\bw^{\top}\bm{\mu}=\rho$ in the Mean-CVaR model \eqref{model:MeanCVaR}. In general, the expected return level $\rho$ is generally set to a positive number. However, this constraint is not always feasible not only in reality but also in theory. The satisfaction of a tight equality constraint for the return level in PO with an ever-changing financial market is difficult. To avoid this problem, we use a regularization term $({\bw}^{\top}\hat{\bm{\mu}}-\rho)^2$ in the objective function instead of using the equality constraint. In this way, the portfolio return ${\bw}^{\top}\hat{\bm{\mu}}$ is balanced between the expected level $\rho$ and the risk metric CVaR. Then the revised $\ell_0$-constrained Mean-CVaR model is given by
\begin{equation}\label{model:mSMCVaR}
\setlength{\abovedisplayskip}{4pt}
\setlength{\belowdisplayskip}{4pt}
\begin{split}
&\min\limits_{\substack{\bw\in\bbRN,\ \tau\in\bbR,\  \bz\in\bbRT}}\left\{\tau+\frac{1}{(1-c)T}\bm{1}^{\top}_{T}\bz+\lambda(\bw^{\top}\bm{\mu}-\rho)^2\right\}\\
&\st\ \bz\geqs-\bR\bw-\tau{\bm1}_T,\ \bz\geqs\bm{0}_{T},\ \bw\in\Delta,\ \|\bw\|_{0}\leqs m.
\end{split}
\end{equation}

For notation simplicity, we define $N_1:=N+1+T$, $N_{2}:=2T+N+2$. We then rewrite model \eqref{model:mSMCVaR} as a more compact form. To this end, we let $\bv:=(\bw^\top,\tau,\bz^\top)^\top\in \bbR^{N_1}$,
\begin{equation}\label{def:h1h2}
\setlength{\abovedisplayskip}{4pt}
\setlength{\belowdisplayskip}{4pt}
\bh_{1}{:=}\bigg(\bm{0}_{N}^\top,1,\frac{1}{(1{-}c)T}\bm{1}_{T}^\top\bigg)^\top,\ \bh_{2}{:=}\left(\bm{\mu}^\top,\bm{0}_{1+T}^\top\right)^\top,
\end{equation}
\begin{equation}
\setlength{\abovedisplayskip}{4pt}
\setlength{\belowdisplayskip}{4pt}
\label{def:tQtI}
\begin{split}
&\tilde{\bQ}:=\left(\begin{array}{ccc}
\bR & \bm{1}_{T} & \bI_{T}\\
\bm{0}_{T\times N} & \bm{0}_{T} & \bI_{T}
\end{array}\right)\in\bbR^{2T\times N_1},\nonumber\\
&\tbI:=\left(\bI_{N}\ \ \bm{0}_{N\times (T+1)}\right)\in\bbR^{N\times N_1},
\end{split}
\end{equation}
and define
\begin{align}\label{eq:defQandq}
\setlength{\abovedisplayskip}{4pt}
\setlength{\belowdisplayskip}{4pt}
&\bQ:=\left(\begin{array}{@{}cc@{}}
\multicolumn{2}{c}{\tilde{\bQ}} \\
\bI_{N} & \bm{0}_{N\times (T+1)} \\
\bm{1}_{N}^\top & \bm{0}_{T+1}^\top \\
\bm{-1}_{N}^\top & \bm{0}_{T+1}^\top \\
\end{array}\right)
\in \bbR^{N_{2} \times N_{1}},\nonumber\\
&\bq:=\left(\begin{array}{@{}c@{}}
\bm{0}_{2T+N} \\
1  \\
-1 \\
\end{array}\right)\in\bbR^{N_2}.
\end{align}
Then model \eqref{model:mSMCVaR} can be rewritten as the following model with respect to (w.r.t.) $\bv$:
\begin{equation}
\setlength{\abovedisplayskip}{4pt}
\setlength{\belowdisplayskip}{4pt}
\label{model:mSMCVaR2}
\begin{split}
&\min\limits_{\bv\in\bbR^{N_1}}\left\{\bh_{1}^{\top} \bv + \lambda (\bh_{2}^{\top}{\bv}-\rho)^2\right\}\\
&\st\ \bQ\bv\geqs\bq,\ \ \|\tilde{\bI}\bv\|_{0}\leqs m.
\end{split}
\end{equation}

Model \eqref{model:mSMCVaR2} comprises two constraints. In fact, it can be reformulated as an equivalent three-term unconstrained minimization model. For this purpose, we define two indicator functions
\begin{equation}
\label{eq:iotaqandiotam}
\begin{split}
&\iota_{\bq}(\bu):=
\begin{cases}
0, & \text{if}\ \bu\geqs\bq;\\
+\infty, & \text{otherwise},
\end{cases}\nonumber\\
&\iota_{m}(\bw):=
\begin{cases}
0, & \text{if}\ \|\bw\|_0\leqs m;\\
+\infty, & \text{otherwise},
\end{cases}
\end{split}
\end{equation}
for $\bu\in\bbR^{N_2}$ and $\bw\in\bbR^{N}$, respectively; and define
\begin{equation}
\label{funcf}
f(\bv):=\bh_{1}^{\top}\bv+\lambda(\bh_{2}^{\top}{\bv}-\rho)^2,\ \ \bv\in\bbR^{N_1}.
\end{equation}
Then model \eqref{model:mSMCVaR2} is equivalent to the following model:
\begin{equation}
\label{model:mSMCVaR3}
\min\limits_{\bv\in\bbR^{N_1}}
\left\{\Phi(\bv):= f(\bv)+
\iota_{\bq}(\bQ\bv)+
\iota_{m}( \tbI \bv)\right\}.
\end{equation}

\begin{figure}[htbp]
\centering
\includegraphics[width=0.9\columnwidth]{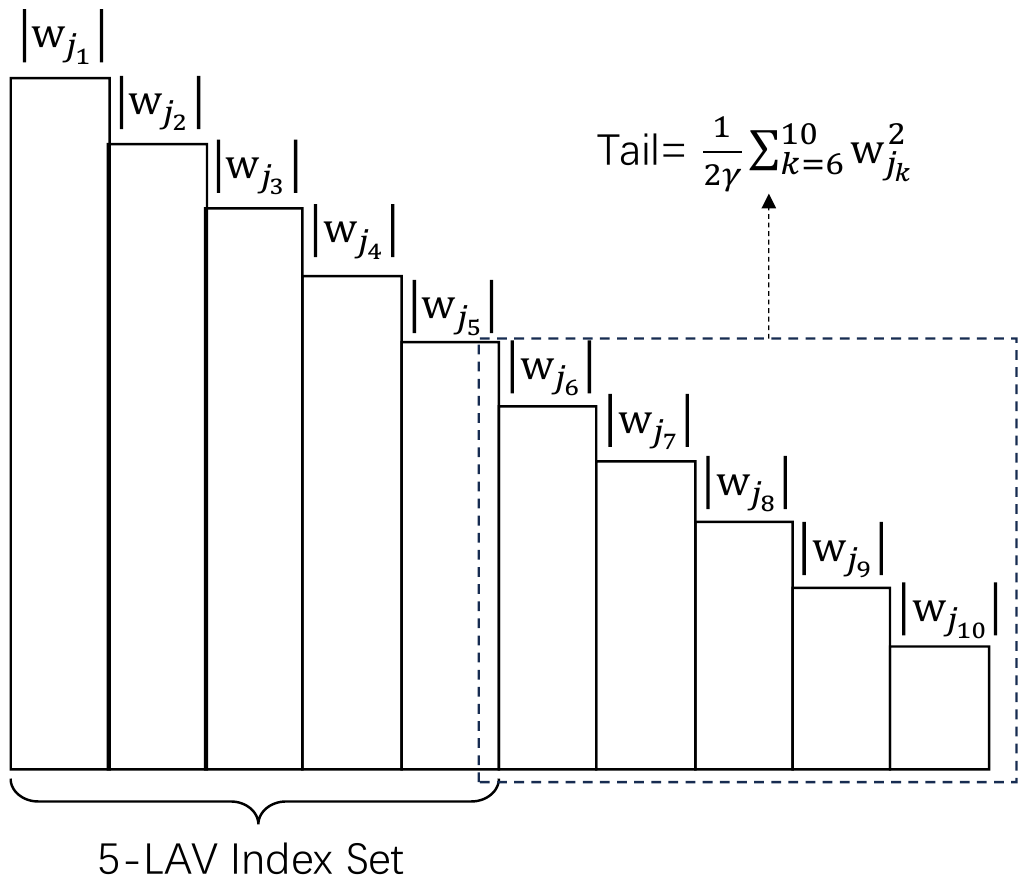}
\caption{Diagram for tailed approximation of $\ell_0$ constraint: $\tilde{\iota}_{m,\gamma}(\bw)$ with $N=10$ and $m=5$. As $\gamma \searrow 0$, $\frac{1}{2\gamma}\sum_{k=6}^{10}w_{j_k}^2 \nearrow +\infty$ and $\tilde{\iota}_{m,\gamma}(\bw)\nearrow \iota_{m}(\bw)$ for any fixed $\bw$ ($\|\bw\|_0>5$).}
\label{fig:tail}
\end{figure}

Directly solving model \eqref{model:mSMCVaR3} is very tough due to the nondifferentiability of $\iota_{\bq}$ and the nonconvexity of $\iota_{m}$. Instead, we introduce a surrogate model of model \eqref{model:mSMCVaR3} and develop an efficient algorithm to solve the surrogate model. To this end, we introduce a function to approximate the exact $m$-sparse function $\iota_m$. For vector $\bw\in\bbRN$, let $\{j_{1},j_{2},\ldots$,$j_{N}\}$ be any rearrangement of $\{1,2,\cdots,N\}$ such that $|w_{j_{1}}|\geqs|w_{j_{2}}|\geqs\cdots\geqs|w_{j_{N}}|$. Then we call $\{ j_{1},j_{2},\cdots,j_{m}\}$ an $m$-largest absolute value ($m$-LAV) index set of $\bw$. Note that the $m$-LAV index sets of $\bw$ may not be unique since the components of $\bw$ may not be distinct. Throughout this paper, for $\bw\in\bbRN$, we denote by $J_{\bw}^m$ an $m$-LAV index set of $\bw$. By smoothing the giant leap of $\iota_{m}$ between $\|\bw\|_0\leqs m$ and $\|\bw\|_0>m$, we construct the following tailed approximation.
\begin{definition}[Tailed Approximation of $\ell_0$ Constraint]\label{def:tailed}
The tailed approximation for the $m$-sparse indicator function is defined by
\begin{equation}
\label{def:iotamgamma}
\tilde{\iota}_{m,\gamma}(\bw)=\begin{cases}
0,&\mbox{if}\ \|\bw\|_0\leqs m;\\
\frac{1}{2\gamma}\sum\limits_{j\in\bbN_N\backslash J_{\bw}^m}w_j^2,&\mbox{if}\ \|\bw\|_0>m.
\end{cases}
\end{equation}
where $\gamma>0$ is the approximation parameter.
\end{definition}
When $\|\bw\|_0>m$, the tail term $\frac{1}{2\gamma}\sum\limits_{j\in\bbN_N\backslash J_{\bw}^m}w_j^2>0$, then $\frac{1}{2\gamma}\sum\limits_{j\in\bbN_N\backslash J_{\bw}^m}w_j^2\nearrow +\infty$ as $\gamma \searrow 0$. Hence for any fixed $\bw$, either $\tilde{\iota}_{m,\gamma}(\bw)= \iota_{m}(\bw)$ ($\|\bw\|_0\leqs m$) or $\tilde{\iota}_{m,\gamma}(\bw)\nearrow \iota_{m}(\bw)$ as $\gamma \searrow 0$ ($\|\bw\|_0>m$). Figure \ref{fig:tail} shows a diagram for $\tilde{\iota}_{m,\gamma}(\bw)$ with $N=10$ and $m=5$.

Now the proposed ASMCVaR model is given by
\begin{equation}
\label{model:ASMCVaR}
\min\limits_{\bv\in\bbR^{N_1}}\left\{\Psi(\bv):=f(\bv)+\iota_{\bq}(\bQ\bv)+\tilde{\iota}_{m,\gamma}(\tbI\bv)\right\}.
\end{equation}
We then establish the relationship between the solutions of model \eqref{model:mSMCVaR3} and the ASMCVaR model. For this purpose, we need the following technical lemma, whose proof is given in Supplementary \ref{proof:soluenviota}.
\begin{lemma}\label{lem:soluenviota}
Let $\bv^{*}$ be a solution of model \eqref{model:ASMCVaR} and $\bw^*:=\tilde{\bI}\bv^*$. Then there exists a constant $\tilde{L}>0$ such that
\begin{equation}
\setlength{\abovedisplayskip}{4pt}
\setlength{\belowdisplayskip}{4pt}
\label{wjcontrbygamma}
w_j^*\leqs\sqrt{2\tilde{L}\gamma},\ \ \mbox{for all}\ j\in\bbN_N\backslash J_{\bw^*}^m.
\end{equation}
\end{lemma}

The following theorem indicates that the optimal value of $f$ in \eqref{model:ASMCVaR} can approximate that in \eqref{model:mSMCVaR3} to arbitrary accuracy by letting $\gamma \searrow 0$.
\begin{theorem}
\label{thm:solutionclose}
Let $\bv^{*}$ and $\bv^{**}$ be solutions of models \eqref{model:ASMCVaR} and \eqref{model:mSMCVaR3}, respectively. Then there exists a constant $\tilde{L}>0$ such that
\begin{equation}
\setlength{\abovedisplayskip}{4pt}
\setlength{\belowdisplayskip}{4pt}
\label{neq:fvoftwomod}
0\leqs f(\bv^{**})-f(\bv^{*})\leqs\sqrt{2\tilde{L}^3\gamma}.
\end{equation}
\end{theorem}
The proof is given in Supplementary \ref{proof:solutionclose}. To solve model \eqref{model:ASMCVaR}, we further rewrite it as an optimization model w.r.t. two variables $\bv$ and $\by$ as follows:
\begin{equation}\label{model:surrogate}
\min\limits_{\bv\in\bbR^{N_1},\ \by\in\bbR^{N}}
\left\{G(\bv,\by):=H(\bv,\by)+\iota_{\bq}(\bQ\bv)+\iota_{m}(\by)\right\},
\end{equation}
where
\begin{equation}
\setlength{\abovedisplayskip}{4pt}
\setlength{\belowdisplayskip}{4pt}
\label{eq:defH}
H(\bv,\by):=f(\bv)+\frac{1}{2\gamma}\left\|\tbI\bv-\by\right\|_{2}^{2},\ \ \gamma>0.
\end{equation}
This model can be directly solved by the proximal alternating linearized minimization algorithm \citep{attouch2013convergence,bolte2014proximal}. We shall then prove the equivalence of the solutions of models \eqref{model:ASMCVaR} and \eqref{model:surrogate}. To this end, we need the following two lemmas. We let
\begin{equation}\label{def:Omegam}
\setlength{\abovedisplayskip}{4pt}
\setlength{\belowdisplayskip}{4pt}
\Omega_m:=\{\bw\in\bbRN:\|\bw\|_0\leqs m\},
\end{equation}
and define the thresholding operator $\mS_{m}$ w.r.t. the $m$ largest components by
\begin{equation}
\setlength{\abovedisplayskip}{4pt}
\setlength{\belowdisplayskip}{4pt}
\label{def:operSm}
[\mS_{m}(\bw)]_j:=\begin{cases}
w_j, &\text{if}\ j\in J_{\bw}^m;\\
0, &\text{if}\ j\in\bbN_N\backslash J_{\bw}^m,
\end{cases}
\end{equation}
for $\bw\in\bbRN$. Note that $\mS_{m}(\bw)$ may not be unique and it depends on a corresponding $m$-LAV index set of $\bw$. Throughout this paper, the equality $\bx=\mS_{m}(\bw)$ represents $\bx\in\mS_{m}(\bw)$.

\begin{lemma}\label{lem:ysSmtIvs}
If a pair $(\bv^{*},\by^{*})$ is a solution of model \eqref{model:surrogate}, then $\by^*=\mS_{m}(\tbI\bv^*)$.
\end{lemma}

\begin{lemma}\label{lem:GvtyleqGvy}
Let $G$ be defined in \eqref{model:surrogate}, $\bv\in\bbR^{N_1}$. If $\tilde{\by}=\mS_m(\tbI\bv)$, then
\begin{equation}\label{neq:GvtyGvy}
\setlength{\abovedisplayskip}{4pt}
\setlength{\belowdisplayskip}{4pt}
G(\bv,\tilde{\by})\leqs G(\bv,\by),\ \forall\by\in\bbRN.
\end{equation}
\end{lemma}

The proof of Lemma \ref{lem:ysSmtIvs} and Lemma \ref{lem:GvtyleqGvy} are given in Supplementary \ref{proof:ysSmtIvs} and \ref{proof:GvtyleqGvy}, respectively. According to these two lemmas, we have the following theorem, whose proof is provided in Supplementary \ref{proof:equivalence}.

\begin{theorem}\label{thm:equivalence}
Let functions $\Psi$ and $G$ be defined by \eqref{model:ASMCVaR} and \eqref{model:surrogate}, respectively. Then
\begin{equation}\label{eq:PsivGvy}
\setlength{\abovedisplayskip}{4pt}
\setlength{\belowdisplayskip}{4pt}
\Psi(\bv)=G(\bv,\by),\ \ \forall\bv\in\bbR^{N_1}, \by=\mS_{m}(\tbI\bv).
\end{equation}
Moreover, a pair $(\bv^{*},\by^{*})$ is a solution of model \eqref{model:surrogate} if and only if $\bv^{*}$ is a solution of model \eqref{model:ASMCVaR} and $\by^*=\mS_{m}(\tbI\bv^*)$.
\end{theorem}

We remark that the original model \eqref{model:mSMCVaR3} presents a three-term single-variable optimization problem (w.r.t. $\bv$), comprising a convex and differentiable term $f$, a convex but nondifferentiable term $\iota_{\bq}\circ{\bQ}$, and a nonconvex and nondifferentiable term $\iota_m\circ\tilde{\bI}$. This model cannot be directly tackled using the PALM algorithm. Resolving such multi-term nondifferentiable optimization models typically necessitates the introduction of proximity operators through subdifferential. When dealing with the summation of two functions $\varphi_1$ and $\varphi_2$, where at least one function is differentiable, the additivity property of subdifferential ($\partial$) holds:
$$
\setlength{\abovedisplayskip}{4pt}
\setlength{\belowdisplayskip}{4pt}
\partial(\varphi_1+\varphi_2)(\bx)=\partial\varphi_1(\bx)+\partial\varphi_2(\bx). $$ If $\varphi_1$ and $\varphi_2$ are both nondifferentiable, this property may not hold. The absence of additivity in subdifferential complicates the design of proximity algorithms to solve the original model. Theoretical proof of convergence becomes exceedingly challenging under such circumstances.

The proposed relaxation technique precisely aims to convert the original model into a three-term dual-variable optimization model (w.r.t. $\bv$ and $\by$) in the form of \eqref{model:surrogate}, which can be solved via the PALM algorithm. The advantage of such a format lies in its ability to decompose the model into two separate optimization problems, each concerning variables $\bv$ and $\by$ respectively. Moreover, due to the differentiability of $H$, the subdifferentials corresponding to these two problems exhibit additivity, thus greatly facilitating the theoretical convergence analysis of the PALM algorithm. Furthermore, we demonstrate, through Theorems \ref{thm:solutionclose} and \ref{thm:equivalence}, the equivalence of solutions in both formats in the sense of approximation, thereby ensuring the effectiveness of utilizing the PALM algorithm to solve the model.

%\newpage
\subsection{Proximal Alternating Linearized Minimization}
In this subsection. We employ the Proximal Alternating Linearized Minimization (PALM) algorithm to solve model \eqref{model:surrogate}. We verify in the following proposition that the partial gradients $\nabla_{\bv}H(\bv,\by)$ and $\nabla_{\by}H(\bv,\by)$ of function $H$ defined by \eqref{eq:defH} are both globally Lipschitz continuous, which is required for the convergence of PALM \citep{bolte2014proximal}.

\begin{proposition}\label{prop:Hproperty}
Let function $H$ be defined by \eqref{eq:defH}. Then there exist positive real numbers $L_1$ and $L_2$ such that
\begin{align*}
\setlength{\abovedisplayskip}{4pt}
\setlength{\belowdisplayskip}{4pt}
\|\nabla_{\bv}H(\bv_1,\by)&-\nabla_{\bv}H(\bv_2,\by)\|_2\leqs L_1\|\bv_1-\bv_2\|_2,\nonumber\\
&\forall\bv_1,\bv_2\in\bbR^{N_1},\ \forall\by\in\bbRN;\\
\|\nabla_{\by}H(\bv,\by_1)&-\nabla_{\by}H(\bv,\by_2)\|_2\leqs L_2\|\by_1-\by_2\|_2,\nonumber\\
&\forall\by_1,\by_2\in\bbR^{N},\ \forall\bv\in\bbR^{N_1}.
\end{align*}
\end{proposition}
The proof is given in Supplementary \ref{proof:Hproperty}. For function $\psi:\bbRn\to(-\infty,+\infty]$, we recall that the proximity operator of $\psi$ at $\bx\in\bbRn$ is defined by
\begin{equation*}
\setlength{\abovedisplayskip}{4pt}
\setlength{\belowdisplayskip}{4pt}
%\label{eq:proxoper}
\prox_{\psi}(\bx):=\argmin\limits_{\bu\in\bbRn}\left\{\frac{1}{2}\|\bu-\bx\|_{2}^{2}+\psi(\bu)\right\}.
\end{equation*}
Then the iterative scheme of PALM to solve model \eqref{model:surrogate} can be given by
\begin{equation*}
\setlength{\abovedisplayskip}{4pt}
\setlength{\belowdisplayskip}{4pt}
%\label{iter:PALM}
\begin{cases}
\bv^{k+1}=\prox_{\iota_{\bq}\circ\bQ}\big(\bv^{k}-\beta_1\nabla_{\bv}H(\bv^{k},\by^{k})\big),\\
\by^{k+1}=\mS_{m}\big(\by^{k}-\beta_2\nabla_{\by}H(\bv^{k+1},\by^{k})\big),
\end{cases}
\end{equation*}
where $\beta_1$ and $\beta_2$ are two introduced step-size parameters. As shown in \citep{bolte2014proximal}, to guarantee the convergence of PALM, it suffices to let $\beta_1\in\left(0,\frac{1}{L_1}\right)$ and $\beta_2\in\left(0,\frac{1}{L_2}\right)$, where $
L_1:=2\lambda\|\bh_{2}\|_2^2+\frac{1}{\gamma}$ and $L_2:=\frac{1}{\gamma}$.

To implement the PALM algorithm, we need the closed form of $\prox_{\iota_{\bq}\circ\bQ}$, which requires solving the following minimization problem:
\begin{equation}
\label{model:proxiotaqQ}
\setlength{\abovedisplayskip}{4pt}
\setlength{\belowdisplayskip}{4pt}
\argmin\limits_{\bu\in\bbR^{N_1}}\left\{\frac{1}{2}\|\bu-\bv\|_{2}^{2}+\iota_{\bq}(\bQ\bu)\right\},
\end{equation}
for a given $\bv\in\bbR^{N_1}$. However, obtaining the closed-form solution of $\prox_{\iota_{\bq}\circ\bQ}$ is very tough, due to the complexity of matrix $\bQ$. Instead, we use the Fixed-Point Proximity Algorithm (FPPA, \citealt{micchelli2011proximity}) to solve model \eqref{model:proxiotaqQ} iteratively. Before this, we give the closed form of operator $\prox_{\iota_{\bq}}$. From the definition of proximity operator, it is easy to see that
\begin{equation*}%\label{eq:iotaqsol}
\setlength{\abovedisplayskip}{4pt}
\setlength{\belowdisplayskip}{4pt}
\prox_{\iota_{\bq}}(\bx){=}\argmin\limits_{\bu\in\left\{\bz\in\bbR^{N_2}:\bz\geqs\bq\right\}}\left\{\frac{1}{2}\|\bu{-}\bx\|_2^{2}\right\}{=}\max\left\{\bx,\bq\right\}.
\end{equation*}
The whole ASMCVaR approach is summarized as Algorithm \ref{alg_ASMCVaR} in Supplementary \ref{supp:algASMCVaR}, with PALM and FPPA in the outer and the inner iterations, respectively.

It is noteworthy that both the proposed relaxation (tailed approximation) technique and the PALM algorithm can be readily extended to various other machine learning problems. To exemplify this, we have introduced two additional modules. Firstly, we apply the proposed relaxation technique and the PALM algorithm to the sparse regression problem, as detailed in Supplementary \ref{app:sparseregress}. Secondly, we have developed a PyTorch module based on this relaxation technique and PALM, facilitating its application to deep learning scenarios. The codes for these two modules are available in the folders \textit{Sparse\_Relaxation\_Test} and \textit{Pytorch\_Demo}, respectively, accessible via the link: \url{https://github.com/linyizun2024/ASMCVaR}.

%\newpage

\section{Experimental Results}
\label{sec:expresult}
We conduct extensive experiments on $6$ real-world financial data sets from Kenneth R. French's Data Library\footnote{\url{http://mba.tuck.dartmouth.edu/pages/faculty/ken.french/data_library.html}}, whose details are provided in Table \ref{tab:datasetsprofiles}. They contain monthly price relative sequences of some portfolios formed by different criteria, and one portfolio here is considered as one asset actually. Both FF25 and FF25EU contain $25$ portfolios, while the former is built on BE/ME (book equity to market equity) and investment from the US market, the latter is built on ME and prior return from the European market. FF32 consists of 32 portfolios built on BE/ME and investment from the US market, while FF49 consists of 49 industry portfolios from the US market. FF100 and FF100MEOP consist of 100 portfolios from the US market. FF100 is built on ME and BE/ME, while FF100MEOP is built on ME and operating profitability, respectively.

\begin{table}[htbp]
\vspace{-5pt}
\setlength{\tabcolsep}{1.2mm}
\centering
\caption{Information of 6 real-world monthly benchmark data sets.}
\scalebox{0.88}{\begin{tabular}{ccccc}
\toprule
Data Set   & Region & Time                   & Months & \# Assets \\ \hline
FF25  & US     & Jul/1971 - May/2023 & 623    & 25     \\
FF25EU      & EU     & Nov/1990 - May/2023 & 391    & 25     \\
FF32      & US     & Jul/1971 - May/2023 & 623    & 32     \\
FF49     & US     & Jul/1971 - May/2023 & 623    & 49    \\
FF100 & US     & Jul/1971 - May/2023 & 623    & 100    \\
FF100MEOP & US     & Jul/1971 - May/2023 & 623    & 100    \\ \bottomrule
\end{tabular}}
\label{tab:datasetsprofiles}
\end{table}

To assess the effectiveness of the proposed ASMCVaR, we take 9 state-of-the-art PO methods as well as 1 baseline into comparison. They are SSMP \citep{sparsepo}, SSPO \citep{SSPO}, SPOLC \citep{SPOLC}, S1, S2, S3 \citep{SSPOl0}, Mean-CVaR \citep{convalueatrisk}, Mean-CVaR-$\ell_1$ \citep{meancvarl1}, MT-CVaR \citep{mtcvar}, and the 1/N strategy (baseline, \citealt{1Nstrategy}). We also try to run the cutting plane approach \citep{cutplanecvar} for Mean-CVaR-$\ell_0$ but the program prompts an error that the models are too large in our experiments. There are 3 slightly different versions in \citep{SSPOl0}: S1 is deterministic but S2 and S3 are randomized. Thus we average the results of 10 times for S2 and S3. The parameters for these competitors are set according to their original papers, while those for ASMCVaR will be set in Section \ref{sec:paraset}.

We adopt the standard moving-window trading scheme in the literature \citep{uPS1,olpsjmlr,egrmvgap} to assess the investing performance of different methods. Assume we observe $\mathscr{T}$ trading periods for $N$ assets in a financial market. Without loss of generality, we can assume the initial wealth $S^{(0)}=1$ and the cumulative wealth (CW) $S^{(t)}$ at time $t$ increases (or decreases) with a factor: $S^{(t)}=S^{(t-1)}\cdot((\mathbf{x}^{(t)})^\top \hat{\bw}^{(t)})$, where $\mathbf{x}^{(t)}:=\br_t+\bm{1}_{N}$ denotes the price relative vector at time $t$ and $\hat{\bw}^{(t)}$ denotes the portfolio for a strategy built on the time window $[t-T:t-1]$. At the beginning where there may be insufficient observations for some strategies, the 1/N strategy can be temporally used. This procedure goes on to the end and forms a  backtest CW sequence $\{S^{(t)}\}_{t=0}^{\mathscr{T}}$ that can be used in the assessment. We set for all the methods $T=60$ as a conventional setting for the time window size in the literature \citep{SSMPsparhed,appmvefficiency}. All the data sets and codes of this section are accessible via the link: \url{https://github.com/linyizun2024/ASMCVaR/tree/main/Codes_for_Experiments_in_Paper}.

\subsection{Parameter Setting and Autonomous Sparsity}
\label{sec:paraset}
The model parameters in \eqref{model:ASMCVaR} are set as follows: the confidence level is set as a conventional one $c=0.99$. The expected return level is empirically set as $\rho=0.02$ based on observations from real-world financial markets. The approximation parameter is set as $\gamma=10^{-5}$, which indicates a tight approximation of $\ell_0$ constraint. The mixing parameter is set as
\begin{equation*}
\begin{split}
\lambda &=  \frac{\sqrt{T}\bh_{1}[(N+2):(N+1+T)]^\top (\bm{1}_{T}/T)}{(\bh_{2}[1:N]^\top (\bm{1}_{N}/N)-\rho)^2}\nonumber\\
&=\frac{1}{(1-c)\sqrt{T}(\bar{r}-\rho)^2},
\end{split}
\end{equation*}
where $\bh_{1}[(N+2):(N+1+T)]$ and $\bh_{2}[1:N]$ denote the last $T$ indices of $\bh_{1}$ and the first $N$ indices of $\bh_{2}$, respectively. $\bar{r}$ is the mean return of the entries of $\bR$. This setting of $\lambda$ is a quotient between the CVaR term and the return term in \eqref{funcf} with equal weights $\bm{1}_{T}/T$ and $\bm{1}_{N}/N$. As for the number of selected assets $m$, we assess $3$ different sparsities: $m_1 = 10$, $m_2 =15$, and $m_3 = 20$.

The algorithm parameters can be conveniently set based on the convergence criteria. For FPPA, we set $\theta=1.99/\|\tbQ\|_2^2$. Its maximum iteration and relative difference tolerance are set as $MaxIter_{1}=200$ and $tol_1=0.001$. For PALM, the learning rates are set as $\beta_1=\frac{0.99}{L_{1}}$ and $\beta_2=\frac{0.99}{L_{2}}$, respectively. Its maximum iteration and relative difference tolerance are set as $MaxIter_{2}=10^{4}$ and $tol_2=10^{-4}$.

We further investigate the sparsity parameter $m$ and find that the support sets of two different sparsities overlap with a large proportion. To be specific, ASMCVaR produces a portfolio $\hat{\bw}^{(t,j)}$ for each sparsity $m_j$ and each trade time $t$. The overlapped proportion can be computed by
\begin{equation*}
\setlength{\abovedisplayskip}{4pt}
\setlength{\belowdisplayskip}{4pt}
p^{(t)}_j:=\frac{|supp(\hat{\bw}^{(t,j)})\cap supp(\hat{\bw}^{(t,j+1)})|}{|supp(\hat{\bw}^{(t,j)})|},\quad j=1,2.
\end{equation*}
The sample means $\bar{p}_j$ and the sample standard deviations (STD) $\hat{\sigma}(p_j)$ of $\{p^{(t)}_j\}_{t=1}^\mathscr{T}$ on the benchmark data sets are presented in Table \ref{tab:sparsityoverlap}. More than $89\%$ assets overlap on average between two different sparsities. We call this property autonomous sparsity since a large proportion of assets can remain in the selected asset pool when the portfolio manager adjusts the pool size. It helps to save transaction costs and reduce the managerial burden.

\vspace{-5pt}
\begin{table}[htbp]
\setlength{\tabcolsep}{1.6mm}
\centering
\caption{Proportion of overlapped assets for ASMCVaR with different sparsities on 6 benchmark data sets.}
\scalebox{0.9}{\begin{tabular}{cccccccccc}
\toprule
   & FF25         & FF25EU              & FF32              & FF49         & FF100  & \begin{tabular}[c]{@{}c@{}}FF100\\ MEOP\end{tabular}  \\
\hline
  $\bar{p}_1$  & 0.9115 & 0.9147 & 0.8919 & 0.9468 & 0.9222 & 0.9018 \\
$\hat{\sigma}(p_1)$	& 0.1182 & 0.1215 & 0.1307 & 0.0832 & 0.1060 & 0.1154 \\
$\bar{p}_2$	& 0.9554 & 0.9553 & 0.9348 & 0.9372 & 0.9543 & 0.9478 \\
$\hat{\sigma}(p_2)$	& 0.0808 & 0.0854 & 0.0958 & 0.0825 & 0.0766 & 0.0724 \\
\bottomrule
\end{tabular}}
\label{tab:sparsityoverlap}
\vspace{-5pt}
\end{table}

\subsection{Cumulative Wealth}
The CW sequence $\{S^{(t)}\}_{t=0}^{\mathscr{T}}$ represents the investing gain of a method along with time. The final CWs of different methods on 6 benchmark data sets are provided in Table \ref{tab:CWs}. ASMCVaR with any of the 3 sparsity settings outperforms all the other competitors on all 6 data sets. For example, the final CW of ASMCVaR is about 10 times those of Mean-CVaR and Mean-CVaR-$\ell_1$ on FF100. Figure \ref{fig:alldatasetCW} in Supplementary \ref{app:alldatasetCW} shows the CW plots of different methods on the benchmark data sets. The CW sequences of ASMCVaR ($m=10$, red color) are steadily over others on most trade times. Besides, ASMCVaR has increasing trends in the long run, indicating that it is growing in wealth during the backtests.

\begin{table*}[htbp]
\setlength{\tabcolsep}{1.5mm}
\footnotesize
\centering
\caption{$\alpha$ factors with p-values of t-tests of different portfolio optimization models on 6 benchmark data sets.}
\label{tab:alphafactor}
\scalebox{0.9}{
\begin{tabular}{cccccccccccccc}
\toprule
\multirow{2}{*}{Method} & \multicolumn{2}{c}{FF25} & \multicolumn{2}{c}{FF25EU} &
\multicolumn{2}{c}{FF32} & \multicolumn{2}{c}{FF49} & \multicolumn{2}{c}{FF100} & \multicolumn{2}{c}{FF100MEOP}\\
\cline{2-13}
&$\alpha$ &p-value &$\alpha$ &p-value &$\alpha$ &p-value &$\alpha$ &p-value &$\alpha$ &p-value  &$\alpha$ &p-value\\
\hline
  1/N  & 0 & 0.4318 & -0.0031 & 1 & 0.0001 & 0.3782 & 0 & 0.5487 & -0.0004 & 0.9427 & -0.0002 & 0.9083 \\
 SSMP   & 0.0002 & 0.4465 & -0.0027 & 0.9992 & -0.0028 & 0.9744 & 0.0017 & 0.1935 & -0.0068 & 0.9999 &  -0.0051 & 0.9983 \\
 SSPO   & -0.0024 & 0.9440 & -0.0098 & 1 & -0.0046 & 0.9994 & -0.0056 & 0.9278 & -0.0100 & 0.9999 &  -0.0071 & 0.9996 \\
  SPOLC  & -0.0032 & 0.9984 & -0.0103 & 1 & -0.0011 & 0.8050 & -0.0037 & 0.9008 & -0.0073 & 1 &  -0.0050 & 0.9977 \\
  S1  & -0.0027 & 0.9701 & -0.0104 & 1 & -0.0049 & 0.9997 & -0.0062 & 0.9463 & -0.0110 & 1 &  -0.0069 & 0.9991 \\
  S2  & -0.0029 & 0.9740 & -0.0101 & 1 & -0.0050 & 0.9996 & -0.0062 & 0.9475 & -0.0106 & 1 & -0.0072 & 0.9994 \\
  S3  & -0.0028 & 0.9738 & -0.0103 & 1 & -0.0049 & 0.9997 & -0.0062 & 0.9465 & -0.0110 & 1 &  -0.0068 & 0.9991 \\
 Mean-CVaR & 0.0009 & 0.1133 & -0.0014 & 0.9216 & 0.0019 & 0.0183 & 0.0025 & 0.0098 & -0.0006 & 0.7222 & 0.0017 & 0.0405 \\
 Mean-CVaR-$\ell_1$ & 0.0015 & 0.0339 & -0.0011 & 0.8947 & 0.0017 & 0.0300 & 0.0027 & 0.0065 & 0 & 0.4987 & 0.0006 & 0.2969 \\
 MT-CVaR & -0.0041 & 0.9984 & -0.0090 & 1 & -0.0067 & 1 & -0.0044 & 0.8979 & -0.0105 & 1 & -0.0050 & 0.9891 \\
 ASMCVaR($m=10$) & \textbf{0.0022} & 0.0006 & \textbf{0.0027} & 0 & \textbf{0.0027} & 0.0002 & \textbf{0.0029} & 0.0049 & \textbf{0.0023} & 0.0038 & \textbf{0.0020} & 0.0060 \\
 ASMCVaR($m=15$) & \textbf{0.0019} & 0.0029 & \textbf{0.0027} & 0 & \textbf{0.0025} & 0.0001 & 0.0024 & 0.0112 & \textbf{0.0025} & 0.0012 & \textbf{0.0017} & 0.0109 \\
 ASMCVaR($m=20$) & \textbf{0.0019} & 0.0017 & \textbf{0.0025} & 0 & \textbf{0.0024} & 0.0001 & 0.0023 & 0.0098 & \textbf{0.0024} & 0.0010 & 0.0016 & 0.0114 \\
\bottomrule
\end{tabular}
}
\label{tb:alphafact}
\end{table*}

\begin{table}[htbp]
\setlength{\tabcolsep}{1.2mm}
\centering
\caption{Final cumulative wealths of different portfolio optimization models on 6 benchmark data sets.}
\scalebox{0.76}{
\begin{tabular}{ccccccccc}
\toprule
Method        & FF25         & FF25EU              & FF32              & FF49         & FF100  & \begin{tabular}[c]{@{}c@{}}FF100\\ MEOP\end{tabular}        \\
\hline
  1/N  & 355.98 & 13.05 & 424.42 & 235.48 & 364.87 & 348.70 \\
 SSMP   & 248.67 & 13.47 & 158.98 & 186.79 & 10.09 & 26.06 \\
 SSPO   & 129.35 & 1.22 & 30.20 & 1.33 & 0.89 & 3.90 \\
  SPOLC  & 57.53 & 0.96 & 169.60 & 5.44 & 2.39 & 8.23 \\
  S1  & 100.76 & 1.08 & 29.47 & 1.09 & 0.54 & 4.31 \\
  S2  & 83.10 & 1.19 & 25.93 & 1.15 & 0.66 & 3.55 \\
  S3  & 96.26 & 1.11 & 29.45 & 1.09 & 0.55 & 4.38 \\
 Mean-CVaR & 333.41 & 18.60 & 331.14 & 150.15 & 70.31 & 306.40 \\
 Mean-CVaR-$\ell_1$ &  440.32 & 18.58 & 332.55 & 178.93 & 102.11 & 163.77 \\
 MT-CVaR & 34.58 & 1.53 & 7.96 & 4.93 & 0.61 & 11.66 \\
 ASMCVaR($m=10$) &	 \textbf{973.23} & \textbf{122.82} & \textbf{1122.76} & \textbf{758.66} & \textbf{966.32} & \textbf{783.47} \\
  ASMCVaR($m=15$) &	 \textbf{827.18} & \textbf{123.13} & \textbf{1303.80} & \textbf{694.83} & \textbf{1126} & \textbf{712.45} \\
  ASMCVaR($m=20$) &	 \textbf{872.03} & \textbf{113.28} & \textbf{1317.01} & \textbf{767.61} & \textbf{1117.20} & \textbf{714.49} \\
\bottomrule
\end{tabular}
}
\label{tab:CWs}
\end{table}

\subsection{$\alpha$ Factor}
According to the Capital Asset Pricing Model (CAPM, \citealt{CAPM}), a portfolio cannot exclude the influence from the underlying financial market. To derive a pure excess return to the market, the $\alpha$ factor is proposed \citep{portalpha0} and becomes a popular score for the relative performance to the market. Let $r_{s}$ and $r_{m}$ denote the
returns for a nontrivial strategy and the Market strategy, respectively. The portfolio return $r_{s}$ can be decomposed into the market return component $r_{m}$ and the $\alpha$ factor:
\begin{equation*}
\setlength{\abovedisplayskip}{4pt}
\setlength{\belowdisplayskip}{4pt}
\bbE\left(r_{s}\right)=\beta \bbE\left(r_{m}\right)+\alpha.
\end{equation*}
$\beta$ and $\alpha$ can be computed by the linear regression:
\begin{equation*}
\setlength{\abovedisplayskip}{4pt}
\setlength{\belowdisplayskip}{4pt}
\hat{\beta}=\frac{\hat{c}\left(r_{s}, r_{m}\right)}{\hat{\sigma}^2\left(r_{m}\right)}, \quad \hat{\alpha}=\bar{r}_{s}-\hat{\beta} \bar{r}_{m} ,
\end{equation*}
where $\hat{c}(\cdot,\cdot)$ and $\hat{\sigma}^2(\cdot)$ are the sample covariance and variance, respectively. $\bar{r}_{s}$ denotes the sample mean, which can be computed by $\bar{r}_{s}=\frac{1}{\mathscr{T}}\sum_{t=1}^\mathscr{T} {r}_{s}^{(t)}$, where ${r}_{s}^{(t)}=(\mathbf{x}^{(t)})^\top \hat{\bw}_s^{(t)}-1$ and $\hat{\bw}_s^{(t)}$ is the portfolio of this nontrivial strategy at time $t$. $\bar{r}_{m}$ can be computed likewise. We adopt the uniform-buy-and-hold strategy \citep{olpsjmlr} as the conventional representative of the market return. To test whether $\alpha>0$ is statistically significant, the right-tailed t-test can also be implemented.

Table \ref{tb:alphafact} shows the $\alpha$ factors with p-values of t-tests of different methods. ASMCVaR($m=10$) achieves the highest $\alpha$s among all the compared methods on all 6 data sets. Moreover, ASMCVaR($m=10$) achieves positive $\alpha$s at the significance level of $1\%$ in all the cases, since all the p-values are smaller than $1\%$. Hence ASMCVaR obtains positive excess returns to the market with high statistical significance.

\subsection{Sharpe Ratio}
Since an investor usually has to take a higher risk to obtain a higher return in financial markets, it is necessary to develop risk-adjusted returns to evaluate this balancing performance. The Sharpe ratio (SR, \citealt{SHARPratio}) is such a metric that takes the ratio between return and risk of a strategy:
\begin{equation}
\setlength{\abovedisplayskip}{4pt}
\setlength{\belowdisplayskip}{4pt}
S R=\frac{\bar{r}_{s}-r_{f}}{\hat{\sigma}\left(r_{s}\right)} ,\notag
\end{equation}
where the risk-free return is set as $r_{f}=0$ in this paper. The (monthly) SRs for different methods on the 6 data sets are shown in Table \ref{tb:sharpe}. ASMCVaR with $m=10$, $m=15$ or $m=20$ achieves the highest SRs among all the competitors on all the data sets. For example, the SRs of ASMCVaR are about $45.7\%$ higher than those of Mean-CVaR and Mean-CVaR-$\ell_1$ on FF25EU.  These results indicate that ASMCVaR performs well in balancing return and risk.

\begin{table}[htbp]
\setlength{\tabcolsep}{1mm}
\footnotesize
\centering
\caption{Sharpe ratios of different portfolio optimization models on 6 benchmark data sets.}
\label{tab:sharpratio}
\scalebox{0.88}{\begin{tabular}{cccccccc}
\toprule
Method & FF25 & FF25EU & FF32 & FF49 & FF100   & \begin{tabular}[c]{@{}c@{}}FF100\\ MEOP\end{tabular}        \\
\hline
  1/N  & 0.2278 & 0.1576 & 0.2236 & 0.2059 & 0.2089 & 0.2077 \\
 SSMP   & 0.1935 & 0.1598 & 0.1536 & 0.1659 & 0.0883 & 0.1085 \\
 SSPO   & 0.1545 & 0.0412 & 0.1182 & 0.0588 & 0.0425 & 0.0685 \\
  SPOLC  & 0.1453 & 0.0316 & 0.1735 & 0.0753 & 0.0563 & 0.0864 \\
  S1  & 0.1498 & 0.0370 & 0.1169 & 0.0559 & 0.0327 & 0.0709 \\
  S2  & 0.1446 & 0.0408 & 0.1137 & 0.0571 & 0.0367 & 0.0670 \\
  S3  & 0.1486 & 0.0381 & 0.1170 & 0.0559 & 0.0331 & 0.0712 \\
 Mean-CVaR & 0.2305 & 0.1803 & 0.2375 & 0.2220 & 0.1738 & 0.2204 \\
 Mean-CVaR-$\ell_1$ & 0.2392 & 0.1862 & 0.2330 & 0.2262 & 0.1815 & 0.1939 \\
 MT-CVaR & 0.1251 & 0.0505 & 0.0844 & 0.0767 & 0.0325 & 0.0911 \\
 ASMCVaR($m=10$) & \textbf{0.2610} & \textbf{0.2758} & \textbf{0.2638} & \textbf{0.2339} & \textbf{0.2446} & \textbf{0.2363} \\
 ASMCVaR($m=15$) & \textbf{0.2541} & \textbf{0.2765} & \textbf{0.2615} & \textbf{0.2286} & \textbf{0.2503} & \textbf{0.2326} \\
 ASMCVaR($m=20$) & \textbf{0.2558} & \textbf{0.2713} & \textbf{0.2603} & \textbf{0.2300} & \textbf{0.2501} & \textbf{0.2322} \\
\bottomrule
\end{tabular}}
\label{tb:sharpe}
\end{table}

\subsection{Transaction Cost}
\label{section:Transaction}
The transaction cost with portfolio change is inevitable in practical portfolio management. We adopt the proportional transaction cost model \citep{UPtc,egrmvgap} to compute the CW at the beginning of the $t$-th trade time:
\begin{align*}
\setlength{\abovedisplayskip}{4pt}
\setlength{\belowdisplayskip}{4pt}
&S^{(\mathscr{T})}{=}S^{(0)} \prod_{t=1}^{\mathscr{T}}\bigg[\big((\mathbf{x}^{(t)})^\top \hat{\bw}^{(t)}\big)\Big(1{-}\frac{\nu}{2} \sum_{i=1}^{N}\big|\hat{w}_{i}^{(t)}{-}\tilde{w}_{i}^{(t-1)}\big|\Big)\bigg],\\
&\tilde{w}_{i}^{(t-1)}=\frac{\hat{w}_{i}^{(t-1)} \cdot \mathrm{x}_{i}^{(t-1)}}{(\mathbf{x}^{(t-1)})^\top \hat{\bw}^{(t-1)}},
\end{align*}
where $\tilde{w}_{i}^{(t-1)}$ denotes the evolved portfolio weight of the $i$-th asset at the end of the $(t-1)$-th trade time, and $\tilde{\bw}^{(0)}$ can be initialized as $\bm{0}_N$. $\nu$ is the bidirectional transaction cost rate, which is examined in the range $[0,0.5\%]$. When the cost rates of buying and selling are the same, changing the evolved portfolio $\tilde{\bw}^{(t-1)}$ to the next portfolio $\tilde{\bw}^{(t)}$ yields a proportional transaction cost of $\frac{\nu}{2} \sum_{i=1}^{N}\big|\hat{w}_{i}^{(t)}-\tilde{w}_{i}^{(t-1)}\big|$.

When $\nu$ changes from $0$ to $0.5\%$, the final CWs for different methods are plotted in Figure \ref{fig:tc} of Supplementary \ref{app:tc}. ASMCVaR outperforms all the other competitors on all the data sets with all $\nu$, showing a good performance against the transaction cost. Note that the 1/N strategy is well known for its little turnover rate mechanism, but ASMCVaR still has an advantage over 1/N, especially with low and medium cost rates. Moreover, if a mutual fund has sufficient capital, then the transaction cost rate $\nu$ can be negotiated to be low. Hence ASMCVaR is effective in practical portfolio management with transaction cost.

\section{Conclusion}
We develop an autonomous sparse mean-CVaR (ASMCVaR) model for portfolio optimization. It aims to optimize a mean-CVaR objective function under the $\ell_{0}$ and the simplex constraints for exact and feasible asset selection. It is an NP-hard problem and the existing combinatorial approach requires a high computational cost. From a markedly different perspective, we propose to convert the $\ell_0$ constraint into an indicator function and handle it with a tailed approximation. By this way, ASMCVaR can approximate the $\ell_0$-constrained mean-CVaR model to arbitrary accuracy. We also develop a proximal alternating linearized minimization algorithm, coupled with a nested fixed-point proximity algorithm (both convergent) to iteratively solve ASMCVaR. It has a good property of autonomous sparsity because a large proportion of assets can remain in the selected asset pool when adjusting the pool size. Experimental results show that ASMCVaR achieves state-of-the-art performance in several common evaluating metrics of investing performance and risk control.

A direct impact of ASMCVaR is to provide a novel tailed approximation approach to the long-standing NP-hard $\ell_0$-constrained problem. It can be conveniently transferred to other related $\ell_0$-constrained problems via similar establishment. Further improvements on ASMCVaR may lie in accelerating the solving algorithm with momentum methods, generalizing the objective function to a broader class of functions, or considering more constraints with practical sense.

\section*{Acknowledgements}
The authors thank the anonymous reviewers for their constructive comments and valuable suggestions in improving this paper. This work was supported in part by National Natural Science Foundation of China under Grant 62176103, in part by Guangdong Basic and Applied Basic Research Foundation under Grant 2021A1515110541, and in part by the Science and Technology Planning Project of Guangzhou under Grants 2024A04J3940 and 2024A04J9896.

% \section*{Impact Statement}
% This paper presents work whose goal is to advance the field of Machine Learning. There are many potential societal consequences of our work, none of which we feel must be specifically highlighted here.

\bibliography{bibfile}

\begin{thebibliography}{34}
\providecommand{\natexlab}[1]{#1}
\providecommand{\url}[1]{\texttt{#1}}
\expandafter\ifx\csname urlstyle\endcsname\relax
  \providecommand{\doi}[1]{doi: #1}\else
  \providecommand{\doi}{doi: \begingroup \urlstyle{rm}\Url}\fi

\bibitem[Ao et~al.(2019)Ao, Li, and Zheng]{appmvefficiency}
Ao, M., Li, Y., and Zheng, X.
\newblock Approaching mean-variance efficiency for large portfolios.
\newblock \emph{The Review of Financial Studies}, 32\penalty0 (7):\penalty0
  2890--2919, 2019.

\bibitem[Artzner et~al.(1999)Artzner, Delbaen, Eber, and Heath]{varsubadd}
Artzner, P., Delbaen, F., Eber, J.-M., and Heath, D.
\newblock Coherent measures of risk.
\newblock \emph{Mathematical Finance}, 9\penalty0 (3):\penalty0 203--228, 1999.

\bibitem[Attouch et~al.(2013)Attouch, Bolte, and
  Svaiter]{attouch2013convergence}
Attouch, H., Bolte, J., and Svaiter, B.~F.
\newblock Convergence of descent methods for semi-algebraic and tame problems:
  proximal algorithms, forward--backward splitting, and regularized
  {Gauss}--{Seidel} methods.
\newblock \emph{Mathematical Programming}, 137\penalty0 (1):\penalty0 91--129,
  2013.

\bibitem[Ban et~al.(2018)Ban, Karoui, and Lim]{mlpo}
Ban, G.-Y., Karoui, N.~E., and Lim, A. E.~B.
\newblock Machine learning and portfolio optimization.
\newblock \emph{Management Science}, 64\penalty0 (3):\penalty0 1136--1154, Mar.
  2018.

\bibitem[Bertsimas \& Cory-Wright(2022)Bertsimas and Cory-Wright]{cutplanePO}
Bertsimas, D. and Cory-Wright, R.
\newblock A scalable algorithm for sparse portfolio selection.
\newblock \emph{INFORMS Journal on Computing}, 34\penalty0 (3):\penalty0
  1489--1511, 2022.

\bibitem[Blum \& Kalai(1999)Blum and Kalai]{UPtc}
Blum, A. and Kalai, A.
\newblock Universal portfolios with and without transaction costs.
\newblock \emph{Machine Learning}, 35\penalty0 (3):\penalty0 193--205, 1999.

\bibitem[Bolte et~al.(2014)Bolte, Sabach, and Teboulle]{bolte2014proximal}
Bolte, J., Sabach, S., and Teboulle, M.
\newblock Proximal alternating linearized minimization for nonconvex and
  nonsmooth problems.
\newblock \emph{Mathematical Programming}, 146\penalty0 (1):\penalty0 459--494,
  2014.

\bibitem[Brodie et~al.(2009)Brodie, Daubechies, Mol, Giannone, and
  Loris]{sparsepo}
Brodie, J., Daubechies, I., Mol, C.~D., Giannone, D., and Loris, I.
\newblock Sparse and stable {M}arkowitz portfolios.
\newblock \emph{Proceedings of the National Academy of Sciences of the United
  States of America}, 106\penalty0 (30):\penalty0 12267--12272, Jul. 2009.

\bibitem[Chao et~al.(2009)Chao, Kavadias, and Gaimon]{POmasci1}
Chao, R.~O., Kavadias, S., and Gaimon, C.
\newblock Revenue driven resource allocation: Funding authority, incentives,
  and new product development portfolio management.
\newblock \emph{Management Science}, 55\penalty0 (9):\penalty0 1556--1569,
  2009.

\bibitem[Cheng \& Gao(2015)Cheng and Gao]{meancvarl1}
Cheng, R. and Gao, J.
\newblock On cardinality constrained mean-cvar portfolio optimization.
\newblock In \emph{The 27th Chinese Control and Decision Conference (2015
  CCDC)}, pp.\  1074--1079, 2015.

\bibitem[Cover(1991)]{uPS1}
Cover, T.~M.
\newblock Universal portfolios.
\newblock \emph{Mathematical Finance}, 1\penalty0 (1):\penalty0 1--29, Jan.
  1991.

\bibitem[DeMiguel et~al.(2009)DeMiguel, Garlappi, and Uppal]{1Nstrategy}
DeMiguel, V., Garlappi, L., and Uppal, R.
\newblock Optimal versus naive diversification: How inefficient is the {1/N}
  portfolio strategy?
\newblock \emph{The Review of Financial Studies}, 22\penalty0 (5):\penalty0
  1915--1953, May 2009.

\bibitem[Dimmock et~al.(2023)Dimmock, Wang, and Yang]{POmasci2}
Dimmock, S.~G., Wang, N., and Yang, J.
\newblock The endowment model and modern portfolio theory.
\newblock \emph{Management Science}, pp.\  1--26, Apr. 2023.
\newblock \doi{10.1287/mnsc.2023.4759}.

\bibitem[Efron et~al.(2004)Efron, Hastie, Johnstone, and Tibshirani]{LARS}
Efron, B., Hastie, T., Johnstone, I., and Tibshirani, R.
\newblock Least angle regression.
\newblock \emph{Annals of Statistics}, 32\penalty0 (2):\penalty0 407--451,
  2004.

\bibitem[Goto \& Xu(2015)Goto and Xu]{SSMPsparhed}
Goto, S. and Xu, Y.
\newblock Improving mean variance optimization through sparse hedging
  restrictions.
\newblock \emph{The Journal of Financial and Quantitative Analysis},
  50\penalty0 (6):\penalty0 1415--1441, 2015.

\bibitem[Gotoh \& Takeda(2011)Gotoh and Takeda]{meancvar}
Gotoh, J.-Y. and Takeda, A.
\newblock On the role of norm constraints in portfolio selection.
\newblock \emph{Computational Management Science}, 8:\penalty0 323--353, 2011.

\bibitem[G\"unl\"uk \& Linderoth(2010)G\"unl\"uk and Linderoth]{MIO}
G\"unl\"uk, O. and Linderoth, J.
\newblock Perspective reformulations of mixed integer nonlinear programs with
  indicator variables.
\newblock \emph{Mathematical Programming}, 124\penalty0 (1):\penalty0 183--205,
  2010.

\bibitem[Jorion(1997)]{valueatrisk}
Jorion, P.
\newblock \emph{Value at Risk: The New Benchmark for Managing Financial Risk}.
\newblock McGraw-Hill, New York, 1997.

\bibitem[Kobayashi et~al.(2021)Kobayashi, Takano, and Nakata]{cutplanecvar}
Kobayashi, K., Takano, Y., and Nakata, K.
\newblock Bilevel cutting-plane algorithm for cardinality-constrained mean-cvar
  portfolio optimization.
\newblock \emph{Journal of Global Optimization}, 81\penalty0 (2):\penalty0
  493--528, 2021.

\bibitem[Lai \& Yang(2023)Lai and Yang]{egrmvgap}
Lai, Z.-R. and Yang, H.
\newblock A survey on gaps between mean-variance approach and exponential
  growth rate approach for portfolio optimization.
\newblock \emph{ACM Computing Surveys}, 55\penalty0 (2):\penalty0 1--36, Mar.
  2023.
\newblock {Article No. 25}.

\bibitem[Lai et~al.(2018)Lai, Yang, Fang, and Wu]{SSPO}
Lai, Z.-R., Yang, P.-Y., Fang, L., and Wu, X.
\newblock Short-term sparse portfolio optimization based on alternating
  direction method of multipliers.
\newblock \emph{Journal of Machine Learning Research}, 19\penalty0
  (63):\penalty0 1--28, Oct. 2018.

\bibitem[Lai et~al.(2020)Lai, Tan, Wu, and Fang]{SPOLC}
Lai, Z.-R., Tan, L., Wu, X., and Fang, L.
\newblock Loss control with rank-one covariance estimate for short-term
  portfolio optimization.
\newblock \emph{Journal of Machine Learning Research}, 21\penalty0
  (97):\penalty0 1--37, Jun. 2020.

\bibitem[Lai et~al.(2022)Lai, Li, Wu, Guan, and Fang]{mtcvar}
Lai, Z.-R., Li, C., Wu, X., Guan, Q., and Fang, L.
\newblock Multitrend conditional value at risk for portfolio optimization.
\newblock \emph{IEEE Transactions on Neural Networks and Learning Systems},
  pp.\  1--14, 2022.
\newblock \doi{10.1109/TNNLS.2022.3183891}.

\bibitem[Li et~al.(2016)Li, Sahoo, and Hoi]{olpsjmlr}
Li, B., Sahoo, D., and Hoi, S.~C.
\newblock {OLPS}: a toolbox for on-line portfolio selection.
\newblock \emph{Journal of Machine Learning Research}, 17\penalty0
  (1):\penalty0 1242--1246, 2016.

\bibitem[Lintner(1965)]{portalpha0}
Lintner, J.
\newblock The valuation of risk assets and the selection of risky investments
  in stock portfolios and capital budgets.
\newblock \emph{Review of Economics and Statistics}, 47\penalty0 (1):\penalty0
  13--37, Feb. 1965.

\bibitem[Liu \& Loewenstein(2013)Liu and Loewenstein]{POmasci3}
Liu, H. and Loewenstein, M.
\newblock Market crashes, correlated illiquidity, and portfolio choice.
\newblock \emph{Management Science}, 59\penalty0 (3):\penalty0 715--732, 2013.

\bibitem[Luo et~al.(2020)Luo, Yu, Xiu, and Wang]{SSPOl0}
Luo, Z., Yu, X., Xiu, N., and Wang, X.
\newblock Closed-form solutions for short-term sparse portfolio optimization.
\newblock \emph{Optimization}, 2020.

\bibitem[Markowitz(1952)]{PS}
Markowitz, H.~M.
\newblock Portfolio selection.
\newblock \emph{Journal of Finance}, 7\penalty0 (1):\penalty0 77--91, Mar.
  1952.

\bibitem[Micchelli et~al.(2011)Micchelli, Shen, and Xu]{micchelli2011proximity}
Micchelli, C.~A., Shen, L., and Xu, Y.
\newblock Proximity algorithms for image models: denoising.
\newblock \emph{Inverse Problems}, 27\penalty0 (4):\penalty0 045009, 2011.

\bibitem[Ogryczak \& Ruszczy\'nski(2002)Ogryczak and Ruszczy\'nski]{ssdorigin}
Ogryczak, W. and Ruszczy\'nski, A.
\newblock Dual stochastic dominance and related mean-risk models.
\newblock \emph{SIAM Journal on Optimization}, 13:\penalty0 60--78, 2002.

\bibitem[Rockafellar \& Uryasev(2000)Rockafellar and Uryasev]{convalueatrisk}
Rockafellar, R.~T. and Uryasev, S.
\newblock Optimization of conditional value-at-risk.
\newblock \emph{Journal of Risk}, 2\penalty0 (3):\penalty0 21--41, 2000.

\bibitem[Sharpe(1964)]{CAPM}
Sharpe, W.~F.
\newblock Capital asset prices: A theory of market equilibrium under conditions
  of risk.
\newblock \emph{Journal of Finance}, 19\penalty0 (3):\penalty0 425--442, Sep.
  1964.

\bibitem[Sharpe(1966)]{SHARPratio}
Sharpe, W.~F.
\newblock Mutual fund performance.
\newblock \emph{Journal of Business}, 39\penalty0 (1):\penalty0 119--138, Jan.
  1966.

\bibitem[Tibshirani(1996)]{LASSO}
Tibshirani, R.
\newblock Regression shrinkage and selection via the lasso.
\newblock \emph{Journal of the Royal Statistical Society}, 58\penalty0
  (1):\penalty0 267--288, 1996.

\end{thebibliography}
\bibliographystyle{icml2024}

\newpage
\appendix
\onecolumn
\section{Supplementary Materials}

\subsection{Proof of Lemma \ref{lem:soluenviota}}
\label{proof:soluenviota}

\begin{proof}
For the case $\|\bw^*\|_0\leqs m$, it is obvious that \eqref{wjcontrbygamma} holds, since $w_j^*=0$ for all $j\in\bbN_N\backslash J_{\bw^*}^m$. We then consider the case $\|\bw^*\|_0>m$. In this case, we construct an approximation $\tilde{\bv}:=(\tilde{\bw}^\top,\tilde{\tau},\tilde{\bz}^\top)^\top\in \bbR^{N_1}$ of $\bv^*:=((\bw^*)^\top,\tau^*,(\bz^*)^\top)^\top$ such that it satisfies the constraints of \eqref{model:mSMCVaR2}, that is, the constraints of \eqref{model:mSMCVaR}. Note that $\bv^*$ satisfies that $\bQ\bv^*\geqs\bq$. We let $\tilde{w}_j=w_j^*$ for $j\in J_{\bw^*}^m\backslash\{j_1\}$, $\tilde{w}_j=0$ for $j\in\bbN_N\backslash J_{\bw^*}^m$ and $\tilde{w}_{j_1}=w_{j_1}^*+\sum\limits_{j\in\bbN_N\backslash J_{\bw^*}^m}w_j^*$. Then $\tilde{\bw}\in\Delta$ and $\|\tilde{\bw}\|_0\leqs m$. For $\tilde{\bz}$ and $\tilde{\tau}$, we define $\delta:=\max\limits_{j\in\bbN_N\backslash J_{\bw^*}^m}\{w_j^*\}$ and $\zeta:=\max\limits_{i\in\bbN_T,j\in\bbN_N}\{|r_{i,j}-r_{i,j_1}|\}$, where $r_{i,j}$ is the $(i,j)$th component of matrix $\bR$ in model \eqref{model:ASMCVaR}, and then let $\tilde{\bz}=\bz^*$, $\tilde{\tau}=\tau^*+\delta\zeta$. It is obvious that $\tilde{\bz}\geqs{\bm0}_T$ since $\bz^*\geqs{\bm0}_T$. To show that $\bQ\tilde{\bv}\geqs\bq$, it remains to be shown that
\begin{equation}\label{neq:ttau1T1}
\tilde{\tau}{\bm1}_T\geqs-\tilde{\bz}-\bR\tilde{\bw}.
\end{equation}
Note that $\bR\bw^*-\bR\tilde{\bw}=\sum\limits_{j\in\bbN_N\backslash J_{\bw^*}^m}w_j^*\left(\br^{(j)}-\br^{(j_1)}\right)$, where $\br^{(j)}$ is the $j$th column of $\bR$. From the definitions of $\delta$ and $\zeta$, we have that
\begin{equation}\label{neq:ttau1T2}
\tilde{\tau}{\bm1}_T\geqs\tau^*{\bm1}_T+\sum_{j\in\bbN_N\backslash J_{\bw^*}^m}w_j^*\left(\br^{(j)}-\br^{(j_1)}\right)=\tau^*{\bm1}_T+\bR\bw^*-\bR\tilde{\bw}.
\end{equation}
Then \eqref{neq:ttau1T1} follows from \eqref{neq:ttau1T2} and the fact $\tau^*{\bm1}_T\geqs-\bz^*-\bR\bw^*$.

Now we know that $\tilde{\iota}_{m,\gamma}(\tilde{\bw})=0$ and $\iota_{\bq}(\bQ\tilde{\bv})=0$. Hence
\begin{equation}\label{eq:Psitvftv}
\Psi(\tilde{\bv})=f(\tilde{\bv}).
\end{equation}
The fact $\bw^*\in\Delta$ gives that $\|\bw^*\|_2\leqs\sqrt{\|\bw^*\|_1}=1$, where $\|\cdot\|_2$ and $\|\cdot\|_1$ denote the $\ell_2$ and $\ell_1$ norms, respectively. Similarly, we also have $\|\tilde{\bw}\|_2\leqs1$. In addition,
$$
\|\bw^*-\tilde{\bw}\|_2^2=\left(\sum_{j\in\bbN_N\backslash J_{\bw^*}^m}w_{j}^*\right)^2+\sum_{j\in\bbN_N\backslash J_{\bw^*}^m}(w_{j}^*)^2\leqs[(N-m)^2+(N-m)]\delta^2
$$
Then it follows from the definition of function $f$ in \eqref{funcf} that
\begin{align}
\notag|f(\bv^*)-f(\tilde{\bv})|&=|(\tau^*-\tilde{\tau})+\lambda({\bm\mu}^\top\bw^*-\rho)^2-\lambda({\bm\mu}^\top\tilde{\bw}-\rho)^2|\\
\notag&\leqs\delta\zeta+\lambda|{\bm\mu}^\top(\bw^*+\tilde{\bw})-2\rho)|\cdot|{\bm\mu}^\top(\bw^*-\tilde{\bw})|\\
\notag&\leqs\delta\zeta+\lambda\left(2\|\bm{\mu}\|_2+2\rho\right)\|\bm{\mu}\|_2\|\bw^*-\tilde{\bw}\|_2\\
\label{neq:fvsftvtLdelta}&\leqs\tilde{L}\delta\leqs\tilde{L},
\end{align}
where
\begin{equation}\label{def:tildeL}
\tilde{L}:=\zeta+2\lambda\sqrt{[(N-m)^2+(N-m)]}\left(\|\bm{\mu}\|_2^2+\rho\|\bm{\mu}\|_2\right).
\end{equation}
This together with \eqref{eq:Psitvftv} and the fact $\Psi(\bv^*)\leqs\Psi(\tilde{\bv})$ yields that
$$
f(\bv^*)+\tilde{L}\geqs f(\tilde{\bv})\geqs\Psi(\bv^*)=f(\bv^*)+\tilde{\iota}_{m,\gamma}(\bw^*)\geqs f(\bv^*)+\frac{1}{2\gamma}\delta^2.
$$
Thus $\delta\leqs\sqrt{2\tilde{L}\gamma}$, which implies the desired result.
\end{proof}

\subsection{Proof of Theorem \ref{thm:solutionclose}}
\label{proof:solutionclose}

\begin{proof}
For the case $\|\tilde{\bI}\bv^*\|_0\leqs m$, since $\tilde{\iota}_{m,\gamma}(\tbI\bv^*)=\iota_{m}(\tbI\bv^{**})=0$ and $\iota_{\bq}(\bQ\bv^{*})=\iota_{\bq}(\bQ\bv^{**})=0$, it is obvious that $f(\bv^{*})=f(\bv^{**})$, and hence \eqref{neq:fvoftwomod} holds. For the case $\|\tilde{\bI}\bv^*\|_0>m$, we define $\tilde{\bv}$ as an approximation of $\bv^*$, following the same definition as in the proof of Lemma \ref{lem:soluenviota}. Then $\bQ\tilde{\bv}\geqs\bq$ and $\|\tilde{\bI}\tilde{\bv}\|_0\leqs m$. Further, $f(\bv^{**})\leqs f(\tilde{\bv})$. By Lemma \ref{lem:soluenviota} and the inequality \eqref{neq:fvsftvtLdelta} in its proof, we know that $|f(\bv^*)-f(\tilde{\bv})|\leqs\tilde{L}\delta\leqs\sqrt{2\tilde{L}^3\gamma}$, where $\tilde{L}$ is defined by \eqref{def:tildeL} and $\delta:=\max\limits_{j\in\bbN_N\backslash J_{\bw^*}^m}\{w_j^*\}$. Thus
\begin{equation}\label{neq:fvssfvsgamma}
f(\bv^{**})\leqs f(\tilde{\bv})\leqs f(\bv^*)+\sqrt{2\tilde{L}^3\gamma}.
\end{equation}
The prerequisite $\bv^{*}$ and $\bv^{**}$ are solutions of models \eqref{model:ASMCVaR} and \eqref{model:mSMCVaR3} gives that $\Psi(\bv^*)\leqs\Psi(\bv^{**})$ and $\Psi(\bv^{**})=f(\bv^{**})$. Hence
\begin{equation}\label{neq:fvsfvss}
f(\bv^*)\leqs\Psi(\bv^*)\leqs f(\bv^{**}).
\end{equation}
Now \eqref{neq:fvoftwomod} can be obtained by \eqref{neq:fvssfvsgamma} and \eqref{neq:fvsfvss} directly. This completes the proof.
\end{proof}

\subsection{Proof of Lemma \ref{lem:ysSmtIvs}}
\label{proof:ysSmtIvs}

\begin{proof}
Since $(\bv^{*},\by^{*})$ is a solution of model \eqref{model:surrogate}, we have that $\|\by^*\|_0\leqs m$ and
\begin{equation}\label{neq:Gvstarystar}
G(\bv^{*},\by^{*})\leqs G(\bv^{*},\by),\ \forall\by\in\bbRN.
\end{equation}
Subtracting $f(\bv^{*})+\iota_{\bq}(\bQ\bv^{*})$ on both sides of \eqref{neq:Gvstarystar} yields that
\begin{equation}\label{neq:iotamystar}
\frac{1}{2\gamma}\|\tbI\bv^{*}-\by^{*}\|_{2}^{2}+\iota_{m}(\by^{*})\leqs\frac{1}{2\gamma}\|\tbI\bv^{*}-\by\|_{2}^{2}+\iota_{m}(\by),\ \forall\by\in\bbRN.
\end{equation}
Let $\bw^*:=\tilde{\bI}\bv^*$. Then \eqref{neq:iotamystar} implies that
\begin{equation}\label{neq:tIvsysy}
\|\bw^{*}-\by^{*}\|_{2}^{2}\leqs\|\bw^{*}-\by\|_{2}^{2},\ \forall\by\in\Omega_m.
\end{equation}
We consider two cases for $\bw^*$. If $\bw^*\in\Omega_{m}$, \eqref{neq:tIvsysy} implies that $\by^*=\bw^*$, that is, $\by^*=\mS_{m}(\bw^*)$.

We then consider the case $\bw^*\notin\Omega_{m}$. In this case, we have $\|\by^*\|_0=m$. Suppose not, then $\|\by^*\|_0<m$. Since $\bw^*\notin\Omega_{m}$, we know that $\bw^*$ has at least $m+1$ nonzero components. Hence there exists some $j'\in\bbN_N$ such that $y_{j'}^*=0$ and $w_{j'}^*\neq0$. By letting $\by\in\Omega_{m}$ such that
\begin{equation}\label{proxiotasetu1}
y_{j'}=w_{j'}^*\ \ \mbox{and}\ \ y_j=y_j^*\ \mbox{for}\ j\in\bbN_N\backslash\{j'\},
\end{equation}
then
\begin{equation}\label{proxiotacontr1}
\|\bw^*-\by^*\|^2-\|\bw^*-\by\|^2=(w_{j'}^*-y_{j'}^*)^2>0,
\end{equation}
which contradicts \eqref{neq:tIvsysy}. Thus $\|\by^*\|_0=m$. Note that the components in $\bw^*$ may not be distinct. We let $J_1,J_2,J_3\subset\bbN_N$ be three index sets such that $|w_{j}^*|>|w_{j_{m}}^*|$ for all $j\in J_1$, $|w_{j}^*|=|w_{j_{m}}^*|$ for all $j\in J_2$, and $|w_{j}^*|<|w_{j_{m}}^*|$ for all $j\in J_3$, respectively.

$(i)$ We first prove that $y_j^*=w_j^*$ for all $j\in J_1$ if $J_1\neq\varnothing$. To this end, we show that $y_j^*\neq 0$ for all $j\in J_1$. Assume that there exists some $j'\in J_1$ such that $y_{j'}^*=0$. Since $\|\by^*\|_0=m$, there must be some $j''\in\bbN_N\backslash J_1$ such that $y_{j''}^*\neq 0$. Note that $|w_{j'}^*|>|w_{j''}^*|$. By letting
$\by\in\Omega_{m}$ such that
\begin{equation}\label{proxiotasetu2}
y_{j'}=w_{j'}^*,\ y_{j''}=0\ \mbox{and}\ y_j=y_j^*\ \mbox{for}\ j\in\bbN_N\backslash\{j',j''\},
\end{equation}
then
\begin{equation}\label{proxiotacontr2}
\|\bw^*-\by^*\|_2^2-\|\bw^*-\by\|_2^2=(w_{j'}^*)^2+(w_{j''}^*-y_{j''}^*)^2-(w_{j''}^*)^2>0,
\end{equation}
which contradicts \eqref{neq:tIvsysy}. Hence $y_j^*\neq 0$ for all $j\in J_1$. Now we are able to prove by contradiction that $y_j^*=w_j^*$ for all $j\in J_1$. Assume that there exists $j'\in J_1$ such that $y_j^*\neq w_{j'}^*$. Since $y_{j'}^*\neq 0$, letting  $\by\in\bbR^N$ be given by \eqref{proxiotasetu1}, then $\by\in\Omega_{m}$ and \eqref{proxiotacontr1} holds, which contradicts \eqref{neq:tIvsysy}. Hence $y_j^*=w_j^*$ for all $j\in J_1$.

$(ii)$ We next prove that $y_j^*=0$ for all $j\in J_3$ if $J_3\neq\varnothing$. Otherwise, there exists $j''\in J_3$ such that $y_{j''}^*\neq0$. Note that the number of elements in $J_1\cup J_2=\bbN_N\backslash J_3$ is greater than or equal to $m$. To guarantee that $\|\by^*\|_0=m$, there must be some $j'\in\bbN_N\backslash J_3$ such that $y_{j'}^*=0$. In this case, we let $\by\in\Omega_{m}$ be given by \eqref{proxiotasetu2}. Note that $|w_{j'}^*|>|w_{j''}^*|$. We have that \eqref{proxiotacontr2} holds, which contradicts \eqref{neq:tIvsysy}. Hence $y_j^*=0$ for all $j\in J_3$.

Note that $J_1\subsetneqq J_{\bw^*}^m$, $J_{\bw^*}^m\backslash J_1\subset J_2$. In addition, $\bw^*\notin\Omega_{m}$ implies that $w_j^*\neq0$ for all $j\in J_1\cup J_2$. Let $m_1\in\{0,1,\ldots,m-1\}$ be the number of elements in $J_1$. The fact $\|\by^*\|_0=m$ together with $(i)$ and $(ii)$ implies that there must be an index set $J_2'\subset J_2$ with $(m-m_1)$ components such that $y_j^*\neq0$ for $j\in J_2'$ and $y_j^*=0$ for $j\in J_2\backslash J_2'$. As shown in $(i)$, we can also prove that $y_j^*=w_j^*$ for $j\in J_2'$. Now by setting $J=J_1\cup J_2'$, then it is an $m$-LAV index set of $\bw$. In addition, we see from \eqref{def:operSm} that $\by^*=\mS_{m}(\tbI\bv^*)$, which completes the proof.
\end{proof}

\subsection{Proof of Lemma \ref{lem:GvtyleqGvy}}
\label{proof:GvtyleqGvy}

\begin{proof}
Let $\bw:=\tbI\bv$. We first prove that
\begin{equation}\label{neq:tywyw}
\|\tilde{\by}-\bw\|_{2}^{2}\leqs\|\by-\bw\|_{2}^{2},\ \forall\by\in\Omega_m.
\end{equation}
If $\bw\in\Omega_{m}$, then it is easy to see from the definition \eqref{def:operSm} of operator $\mS_m$ that $\tilde{\by}=\bw$ and hence \eqref{neq:tywyw} holds. For the case $\bw\notin\Omega_{m}$, since $\tilde{\by}=\mS_m(\tbI\bv)$, we know that $\tilde{\by}\in\Omega_{m}$ and
\begin{equation}\label{eq:tywsumwj}
\|\tilde{\by}-\bw\|_2^2=\sum_{j\in\bbN_N\backslash J_{\bw}^m}w_j^2.
\end{equation}
For any $\by\in\Omega_{m}$, there exists an index set $J_{\by}\subset\bbN_N$ with $m$ elements such that $y_j=0$ for all $j\in\bbN_N\backslash J_{\by}$. Since $J_{\bw}^m$ contains the $m$ largest components of $\bw$ in absolute value, we have that $\sum\limits_{j\in J_{\bw}^m}w_j^2\geqs\sum\limits_{j\in J_{\by}}w_j^2$, which together with \eqref{eq:tywsumwj} implies that
$$
\|\by-\bw\|_2^2\geqs\sum_{j\in\bbN_N\backslash J_{\by}}w_j^2\geqs\sum_{j\in\bbN_N\backslash J_{\bw}^m}w_j^2=\|\tilde{\by}-\bw\|_2^2.
$$
This completes the proof of \eqref{neq:tywyw}. Now by \eqref{neq:tywyw}, we have
\begin{equation}\label{neq:iotamty}
\frac{1}{2\gamma}\|\tilde{\by}-\tbI\bv\|_{2}^{2}+\iota_{m}(\tilde{\by})\leqs\frac{1}{2\gamma}\|\by-\tbI\bv\|_{2}^{2}+\iota_{m}(\by),\ \forall\by\in\bbRN.
\end{equation}
Adding $f(\bv)+\iota_{\bq}(\bQ\bv)$ on both sides of \eqref{neq:iotamty} yields \eqref{neq:GvtyGvy} immediately.
\end{proof}

\subsection{Proof of Theorem \ref{thm:equivalence}}
\label{proof:equivalence}

\begin{proof}
Let $\bw:=\tilde{\bI}\bv$. The definition of $\tilde{\iota}_{m,\gamma}$ in \eqref{def:iotamgamma} implies that $\tilde{\iota}_{m,\gamma}(\bw)=\frac{1}{2\gamma}\|\bw-\mS_{m}(\bw)\|_2^2$. Then for $\by=\mS_{m}(\bw)$, we know that $\iota_m(\by)=0$ and $\tilde{\iota}_{m,\gamma}(\bw)=\frac{1}{2\gamma}\|\bw-\by\|_2^2$. Hence $\Psi(\bv)=G(\bv,\by)$.

Suppose that $(\bv^{*},\by^{*})$ is a solution of model \eqref{model:surrogate}. From Lemma \ref{lem:ysSmtIvs}, we know that $\by^*=\mS_{m}(\tbI\bv^*)$, and hence $\Psi(\bv^*)=G(\bv^*,\by^*)$. If $\bv^{*}$ is not a solution of model \eqref{model:ASMCVaR}, then there exists $\tilde{\bv}\in\bbR^{N_1}$ such that $\Psi(\tilde{\bv})<\Psi(\bv^{*})$. In this case, we have that $G(\tilde{\bv},\tilde{\by})=\Psi(\tilde{\bv})<\Psi(\bv^{*})=G(\bv^{*},\by^*)$, for any $\tilde{\by}=\mS_{m}(\tbI\tilde{\bv})$, which contradicts the assumption that $(\bv^{*},\by^{*})$ is a solution of model \eqref{model:surrogate}. This proves that $\bv^{*}$ is a solution of model \eqref{model:ASMCVaR}.

Conversely, suppose that $\bv^{*}$ is a solution of model \eqref{model:ASMCVaR} and $\by^{*}=\mS_m(\tbI{\bv^{*}})$. We confirm that $(\bv^{*},\by^{*})$ is a solution of model \eqref{model:surrogate}. If it is not true, according to Lemma \ref{lem:GvtyleqGvy}, we can find vectors $\bar{\bv}\in\bbR^{N_1}$ and $\bar{\by}=\mS_m(\tbI{\bar{\bv}})$ such that $G(\bar{\bv},\bar{\by})<G(\bv^{*},\by^{*})$. We see from \eqref{eq:PsivGvy} that $\Psi(\bar{\bv})=G(\bar{\bv},\bar{\by})$ and $\Psi(\bv^*)=G(\bv^*,\by^*)$. Hence $\Psi(\bar{\bv})<\Psi(\bv^*)$, which contradicts that $\bv^{*}$ is a solution of model \eqref{model:ASMCVaR}. This completes the proof.
\end{proof}

\subsection{Proof of Proposition \ref{prop:Hproperty}}
\label{proof:Hproperty}

\begin{proof}
Compute that
\begin{align*}
%\label{eq:Lipv}
&\nabla_{\bv} H(\bv,\by)=\bh_{1}+2\lambda\left(\bh_{2}\bh_{2}^{\top}\bv-\rho\bh_{2}\right) +\frac{1}{\gamma}\left(\tbI^\top\tbI\bv-\tbI^\top\by\right),\\
%\label{eq:Lipy}
&\nabla_{\by}H(\bv,\by)=\frac{1}{\gamma}(\by-\tbI\bv).
\end{align*}
Let
\begin{equation*}%\label{def:LipconstL1L2}
L_1:=2\lambda\|\bh_{2}\|_2^2+\frac{1}{\gamma}\ \ \mbox{and}\ \ L_2:=\frac{1}{\gamma}.
\end{equation*}
Then for all $\bv_1,\bv_2\in\bbR^{N_1}$ and $\by\in\bbRN$,
\begin{align*}
\|\nabla_{\bv}H(\bv_1,\by)-\nabla_{\bv}H(\bv_2,\by)\|_2&=\left\|\left(2\lambda\bh_{2}\bh_{2}^{\top}+\frac{1}{\gamma}\tbI^{\top}\tbI\right)(\bv_{1}-\bv_{2})\right\|_2\\
&\leqs\left\|\left(2\lambda\bh_{2}\bh_{2}^{\top}+\frac{1}{\gamma}\tbI^{\top}\tbI\right)\right\|_2\cdot\|\bv_{1}-\bv_{2}\|_2\\
&\leqs\left(2\lambda\|\bh_{2}\|_2^2+\frac{1}{\gamma}\right)\|\bv_{1}-\bv_{2}\|_2\\
&=L_{1}\|\bv_{1}-\bv_{2}\|_2.
\end{align*}
For all $\by_1,\by_2\in\bbR^{N}$ and $\bv\in\bbR^{N_1}$,
$$
\|\nabla_{\by}H(\bv,\by_1)-\nabla_{\by}H(\bv,\by_2)\|_2=\frac{1}{\gamma}\|\by_1-\by_2\|_2=L_2\|\by_1-\by_2\|_2,
$$
which completes the proof.
\end{proof}

\subsection{Algorithm for ASMCVaR}\label{supp:algASMCVaR}

\begin{algorithm}
\caption{ASMCVaR}\label{alg_ASMCVaR}
\begin{algorithmic}
\REQUIRE Given the sample asset price relative matrix $\bX\in\bbR^{T\times N}$; the parameters $m$ and $\lambda$; the maximum iteration numbers $(MaxIter_1,MaxIter_2)$ and tolerances $(tol_1,tol_2)$. Set the approximation parameter $\gamma$, the expected return level $\rho$ and the confidence level $c$.\\
\textbf{Preparation}: Compute the return matrix $\bR=\bX-\bm{1}_{T\times N}$, the sample mean vector $\hat{\bm{\mu}}=\frac{1}{T}\bR^{\top}\bm{1}_{T}$. Let matrix $\bQ$ and vector $\bq$ be given by \eqref{eq:defQandq}, $\tbI$ be given by \eqref{def:tQtI}, $\bh_{1}$ and $\bh_{2}$ be given by \eqref{def:h1h2}. Compute the Lipschitz constants $L_1=2\lambda\|\bh_{2}\|_2^2+\frac{1}{\gamma}$ and $L_2=\frac{1}{\gamma}$.\\
\textbf{Initialization}: Given the initial vectors by $\bv^{0}=\left(\frac{1}{N}\bm{1}_{N}^\top,\bm{0}_{1+T}^\top\right)^\top$ and $\by^{0}=\frac{1}{N}\bm{1}_{N}$. Set $\beta_1=\frac{0.99}{L_{1}}$ and $\beta_2=\frac{0.99}{L_{2}}$ in PALM. Set $\theta=\frac{1.99}{\|\bQ\|_2^{2}}$ in FPPA.\\
\textbf{repeat}\\
$\bp^k{=}\bv^{k}{-}\beta_1\left[\bh_{1}{+}2\lambda\left(\bh_{2}\bh_{2}^{\top}\bv^k{-}\rho\bh_{2}\right) {+}\frac{1}{\gamma}\left(\tbI^\top\tbI\bv^k{-}\tbI^\top\by^k\right)\right]$

\hspace{2em}\textbf{repeat}\\
1.\hspace{3em}$\bx^{l}=\bQ\bp^k+\bz^{l}-\theta\bQ\bQ^\top\bz^{l}$\\
2.\hspace{3em}$\bz^{l+1}=\bx^{l}-\max\left\{\bx^{l},\bq\right\}$\\
3.\hspace{3em}$l=l+1$

\hspace{2em}\textbf{until} $\frac{\|\bz^{l}-\bz^{l-1}\|_2}{\|\bz^{l-1}\|_2}\leqs tol_1$ or $l>MaxIter_1$.\\
4.\hspace{1.2em}$\bv^{k+1}=\bp^k-\theta\bQ^\top\bz^{l}$\\
5.\hspace{1.2em}$\by^{k+1}=\mS_{m}\left(\by^{k}-\beta_2\left(\frac{1}{\gamma}(\by^k-\tbI\bv^{k+1})\right)\right)$\\
6.\hspace{1.2em}$k=k+1$\\
\textbf{until} $\frac{\|\bv^{k}-\bv^{k-1}\|_{2}}{\|\bv^{k-1}\|_{2}}\leqs tol_{2}$ or $k>MaxIter_{2}$.\\
\ENSURE Portfolio $\bw^{*}=\bv_{1:N}^k$.
\end{algorithmic}
\end{algorithm}

\subsection{PALM for Relaxed Sparse Regression}\label{app:sparseregress}

To see the performance of the proposed relaxation technique and the PALM algorithm for the sparse regression problem, we consider the following two models:
$$
\min_{\bm{\beta}}\left\{\frac{1}{2}\|{\bX}{\bm{\beta}}-\by\|_2^2+\iota_m(\bm{\beta})\right\},
$$
$$
\min_{\bm{\beta},\bm{\eta}}\left\{\frac{1}{2}\|\bX\bm{\beta}-\by\|_2^2+\frac{1}{2\gamma}\|\bm{\beta}-\bm{\eta}\|_2^2+\iota_m(\bm{\eta})\right\}.
$$
They correspond to the original model \eqref{model:mSMCVaR3} and the relaxed dual-variable optimization model \eqref{model:surrogate}, respectively. We conducted two sets of experiments to validate two aspects:
\begin{itemize}[topsep=0pt, itemsep=0pt, parsep=0pt]
\item[1.] The equivalence of the above relaxation model and original model, in the sense of approximation.
\item[2.] The effectiveness of the PALM algorithm for solving the relaxation model.
\end{itemize}

In these experiments, data sets comprising 50 samples are synthetically generated according to the model $y=\bx^{\top}\bm{\beta}^*+\epsilon$, where $\bm{\beta}^*=({\bm1}_3,{\bm0}_7)^\top\in\bbR^{10}$ represents the true coefficient vector and $\epsilon$ denotes the random error following a normal distribution $\mN(0,1)$. The input vector $\bx\in\bbR^{10}$ is generated by a normal distribution with mean zero and covariance matrix $\bm{\Sigma}\in\bbR^{10\times 10}$, where $\Sigma_{ij}:=0.5^{|i-j|}$. We set the sparsity level $m$ as 3. The direct exhaustive approach is utilized to enumerate all possible support set configurations, totaling $C_{10}^{3}=120$ cases. By comparing the optimal solutions corresponding to these 120 cases, we can ascertain the exact globally optimal solutions of the aforementioned two models. The iterative scheme of the PALM algorithm for solving the above relaxed model can be given by
\begin{equation*}
\begin{cases}
\bm{\beta}^{k+1}=\bm{\beta}^{k}-\alpha_1\left[\bX^\top(\bX\bm{\beta}^k-\by)+\frac{1}{\gamma}(\bm{\beta}^k-\bm{\eta}^k)\right],\\
\bm{\eta}^{k+1}=\mS_m\left(\left(1-\frac{\alpha_2}{\gamma}\right)\bm{\eta}^k+\frac{\alpha_2}{\gamma}\bm{\beta}^{k+1}\right),
\end{cases}
\end{equation*}
where $\mS_m$ is defined by \eqref{def:operSm}. Parameters $\alpha_1$ and $\alpha_2$ are set as $0.99\cdot\|\bX^\top\bX+\frac{1}{\gamma}\bI\|_2^{-1}$ and $0.99\gamma$, respectively. To ensure a better convergence, we conduct $5\times 10^6$ iterations of PALM. The initial vectors for iteration are casually chosen as $\bm{\beta}^0=\bm{\eta}^0=\bm{0}_{10}$.

To demonstrate the robustness of our findings, we repeated each of the aforementioned two sets of experiments for 10 times, with different $\bX$ and $\by$ in each run. Through the first set of experiments, we observed that for sufficiently small $\gamma$ (e.g., $10^{-7}$), the solutions of the original and relaxed models are identical. This suggests that as $\gamma$ tends to zero, the optimal solution of the relaxed model indeed converges to that of the original model. In the subsequent set of experiments, we observed that across the 10 trials, the solutions obtained by PALM consistently approached the optimal solutions of the relaxed model to a remarkable degree. This reflects that in certain scenarios, PALM may be capable of directly converging to the globally optimal solution of the original nonconvex optimization model, even though the initializations may not necessarily be very close to the true solution.

\begin{figure*}[htbp]
\centering
\subfloat[FF25]{
\includegraphics[width=0.45\textwidth]{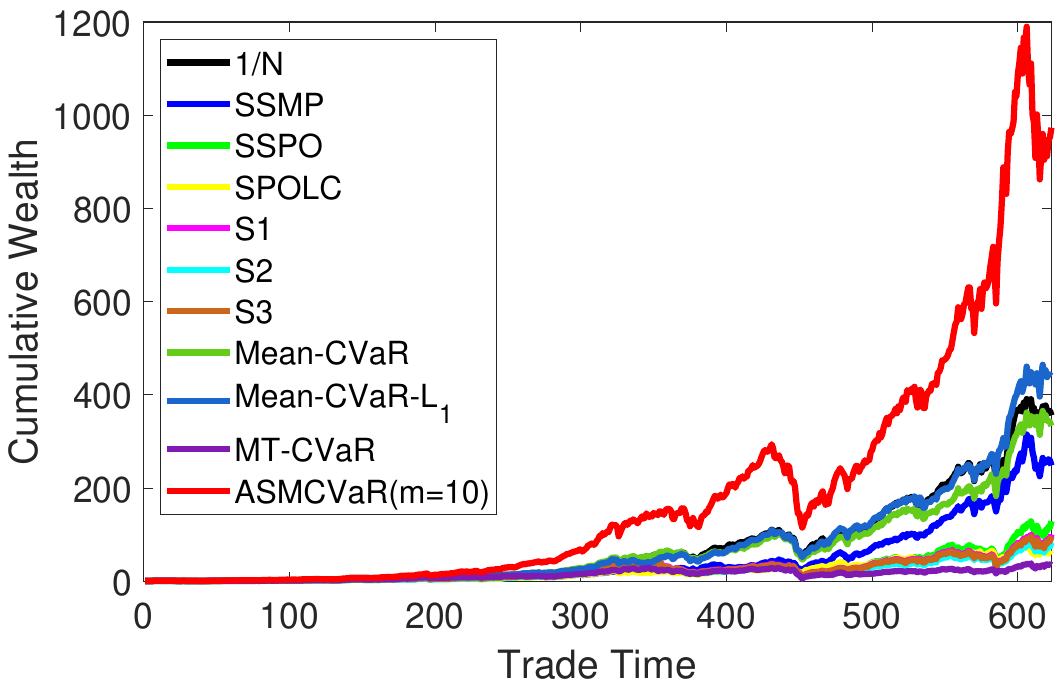}}
\hspace{4pt}\subfloat[FF25EU]{
\includegraphics[width=0.45\linewidth]{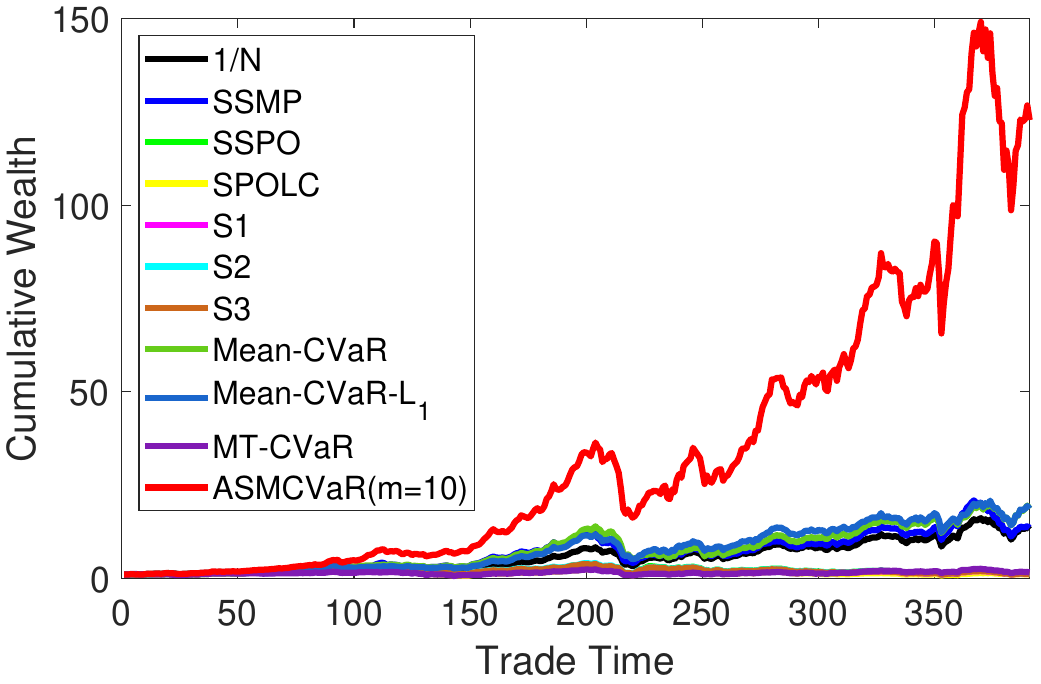}}\\
\subfloat[FF32]{
\includegraphics[width=0.45\linewidth]{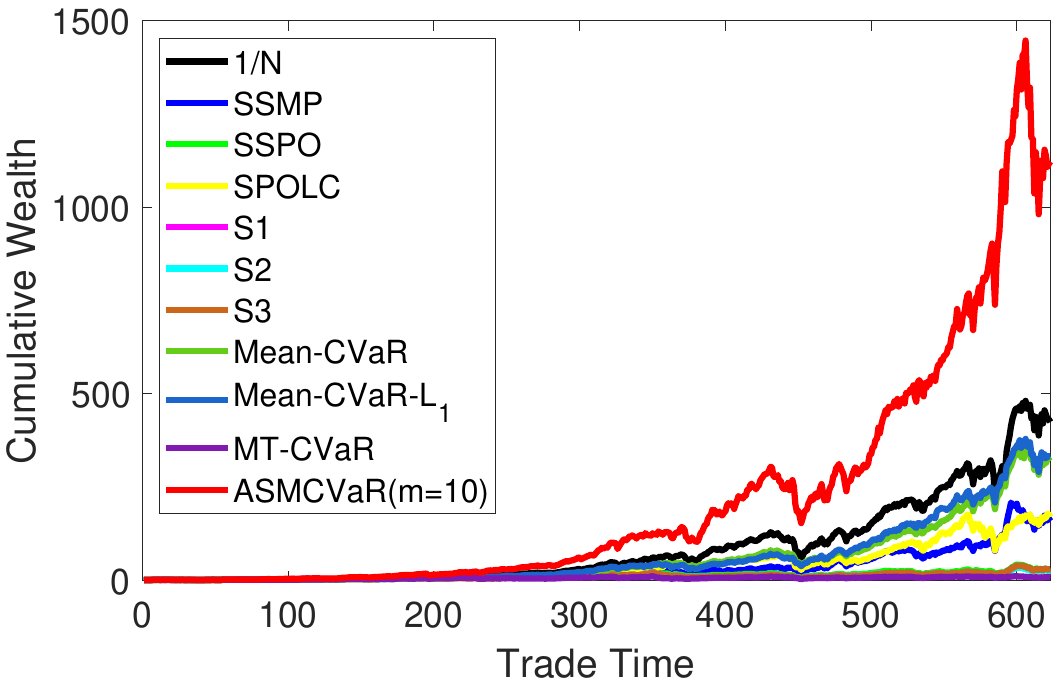}}
\subfloat[FF49]{
\hspace{4pt}\includegraphics[width=0.45\linewidth]{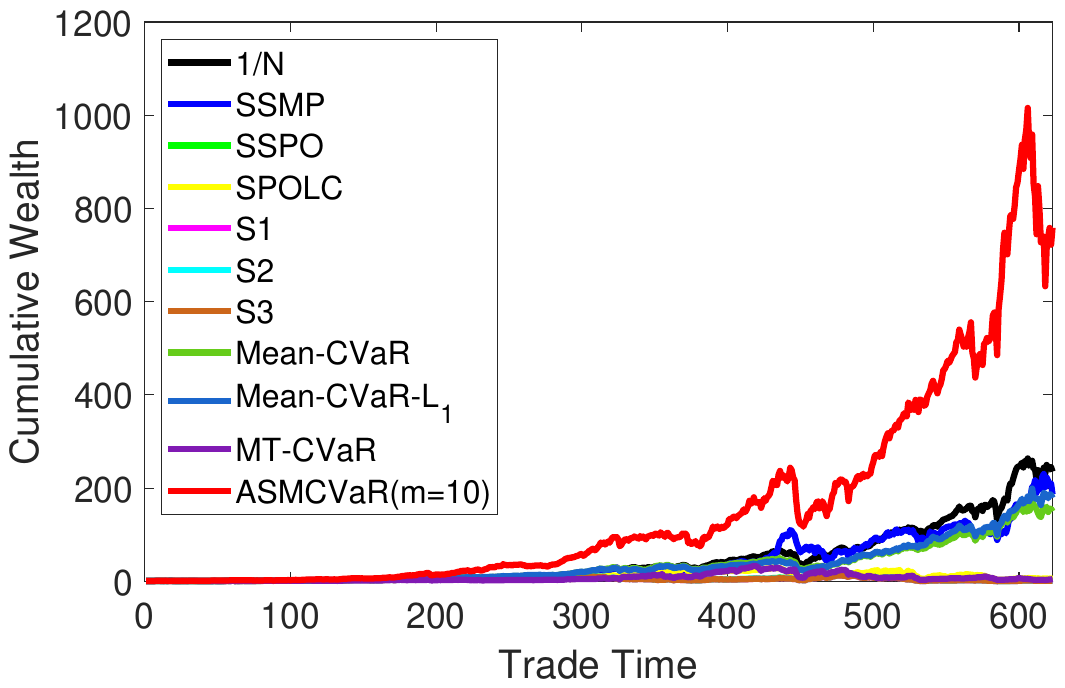}}\\
\subfloat[FF100]{
\includegraphics[width=0.45\linewidth]{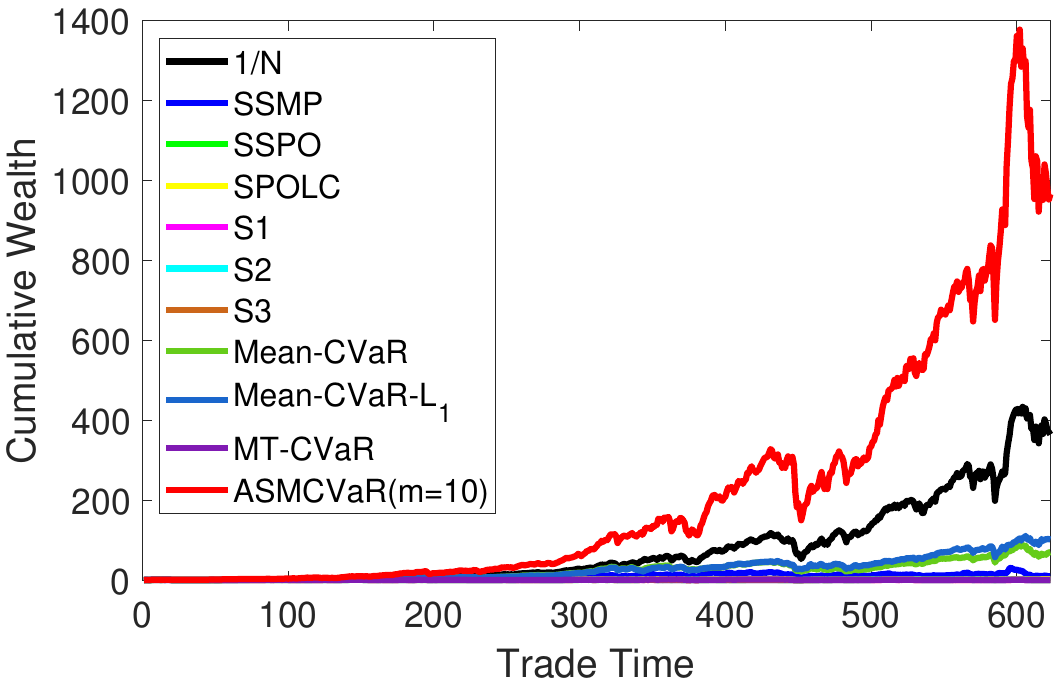}}
\subfloat[FF100MEOP]{
\hspace{4pt}\includegraphics[width=0.45\linewidth]{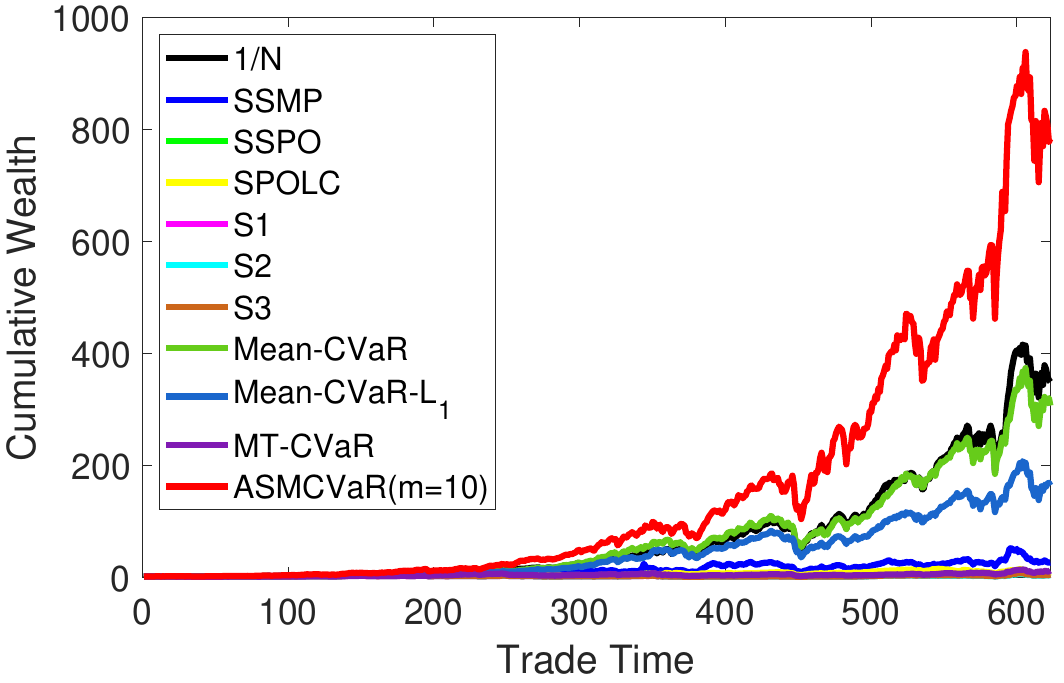}}
\caption{Cumulative wealths of different portfolio optimization models along with trade time on 6 benchmark data sets.}
\label{fig:alldatasetCW}
\end{figure*}

\subsection{Cumulative Wealth Plots (Figure \ref{fig:alldatasetCW})}
\label{app:alldatasetCW}

\begin{figure*}[htbp]
\centering
\subfloat[FF25]{
\centering
\includegraphics[width=0.42
\textwidth]{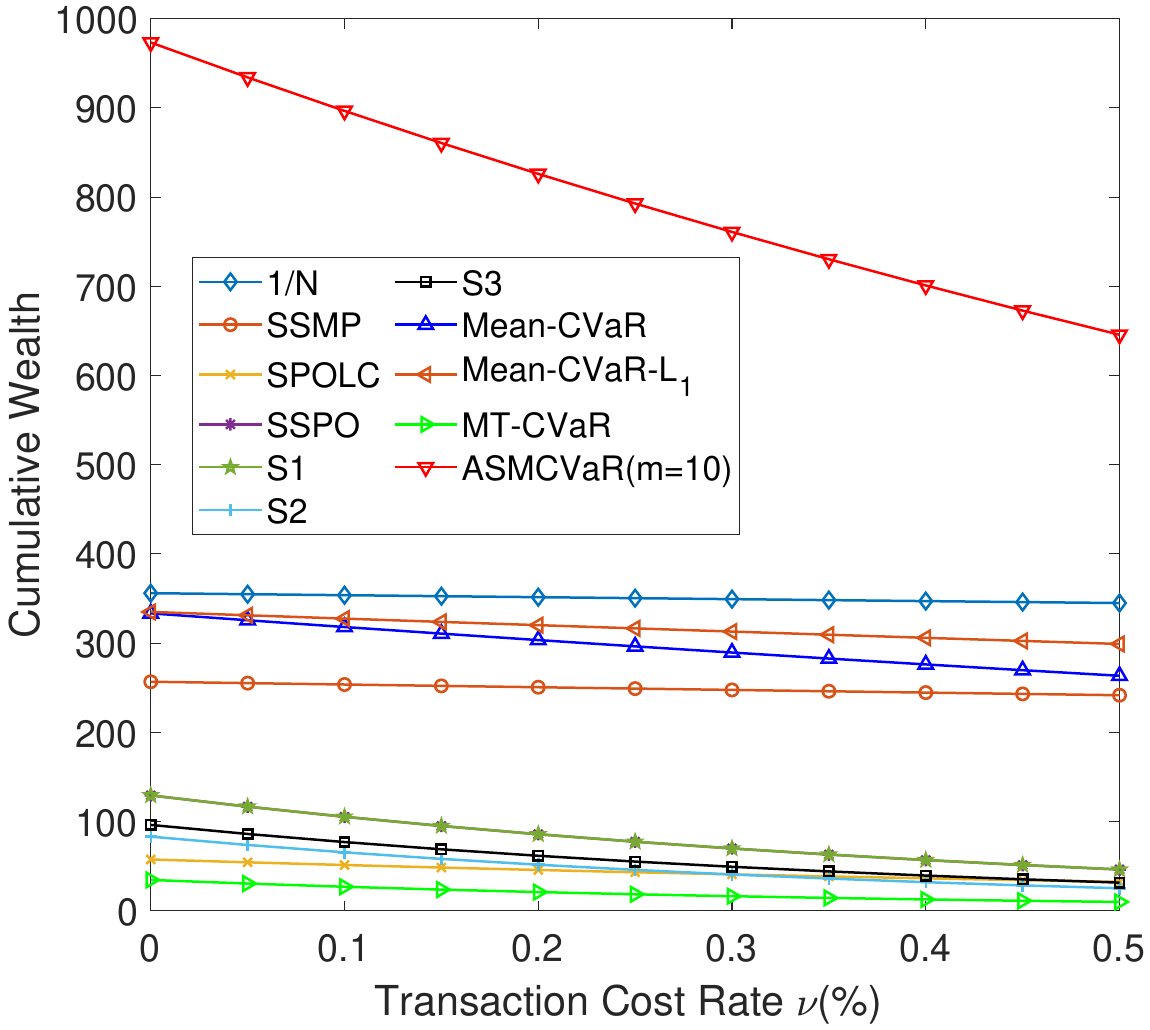}}
\subfloat[FF25EU]{
\centering
\includegraphics[width=0.42\textwidth]{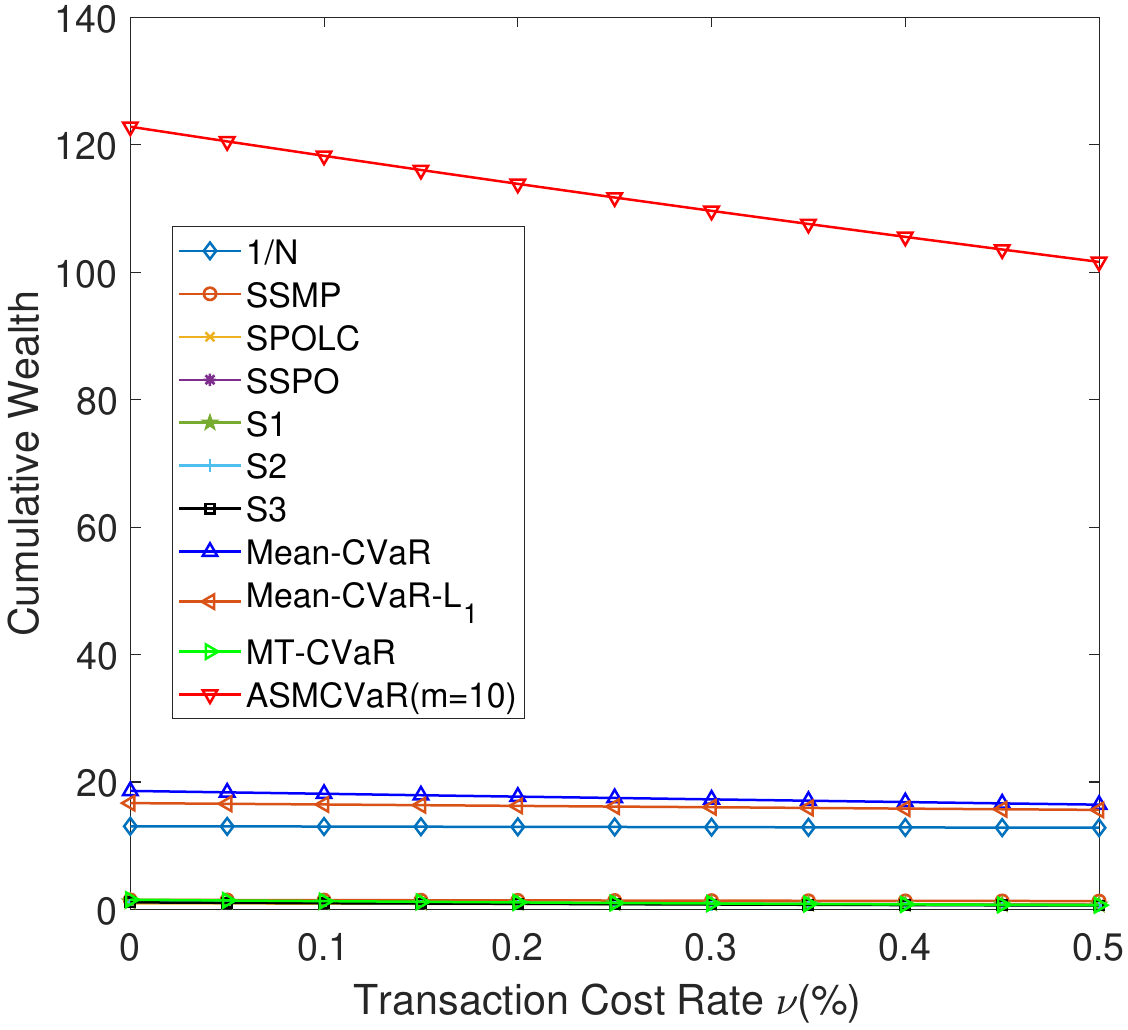}}\\\vspace{-0.5em}
\subfloat[FF32]{
\centering
\includegraphics[width=0.42\textwidth]{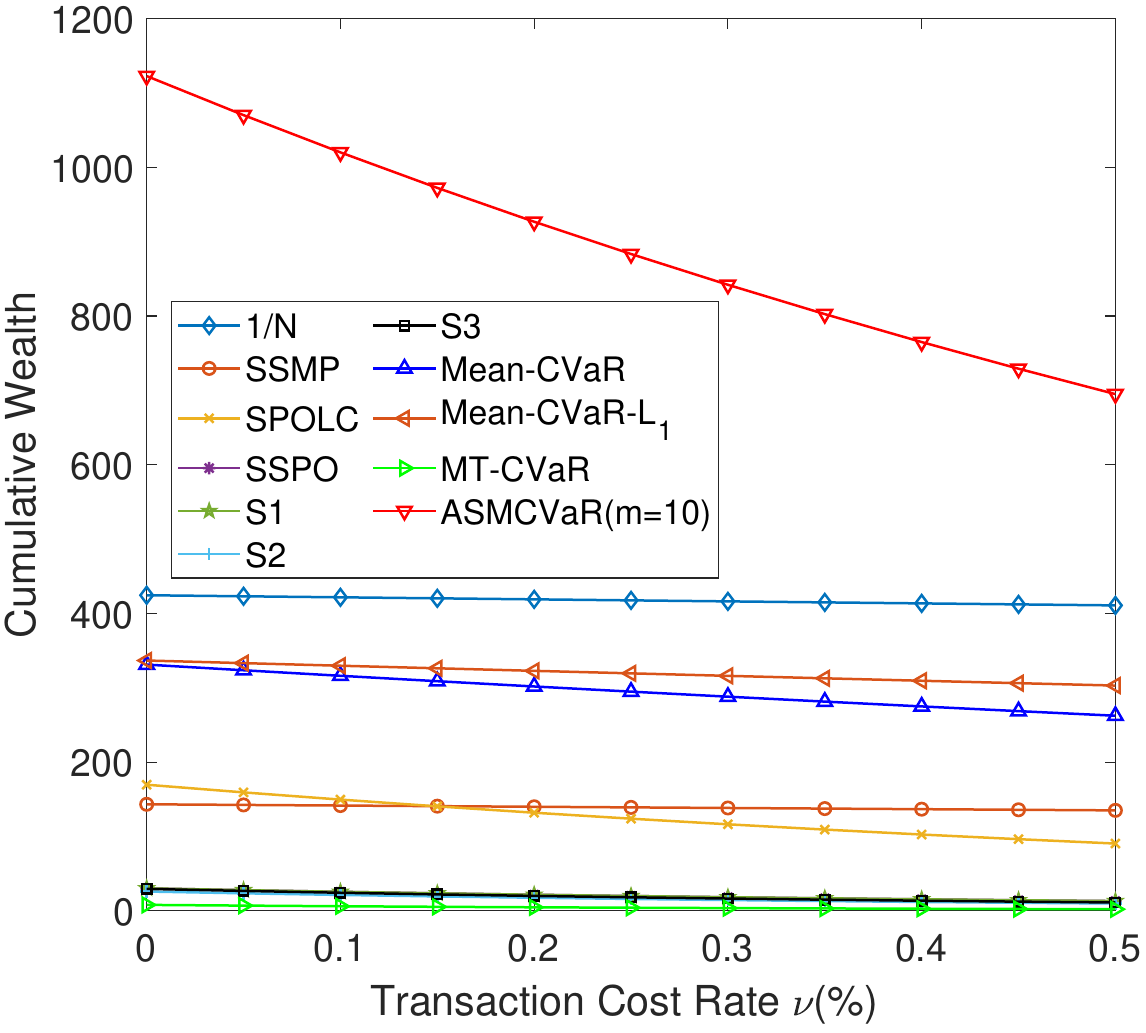}}
\subfloat[FF49]{
\centering
\includegraphics[width=0.42\textwidth]{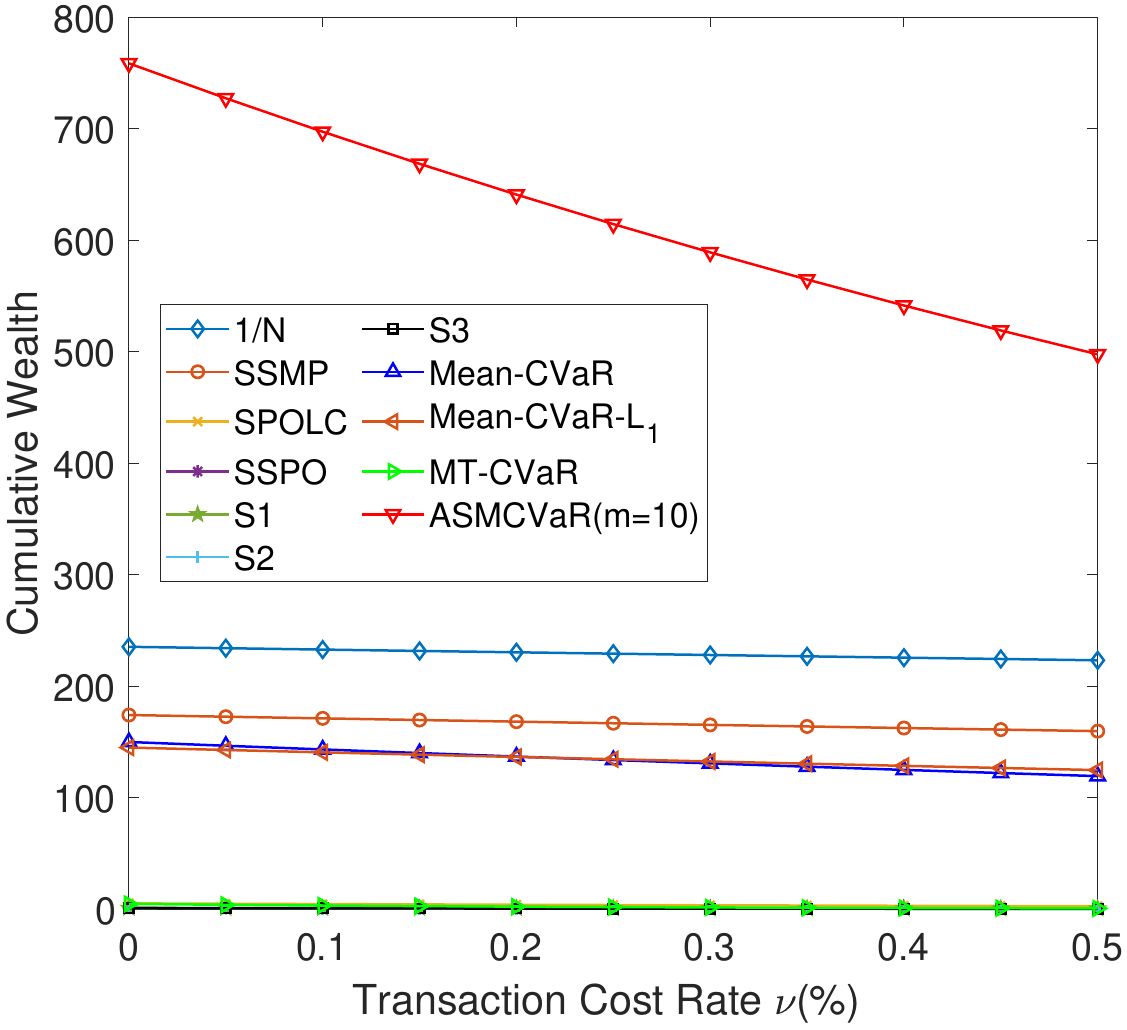}}\\\vspace{-0.5em}
\subfloat[FF100]{
\centering
\includegraphics[width=0.42\textwidth]{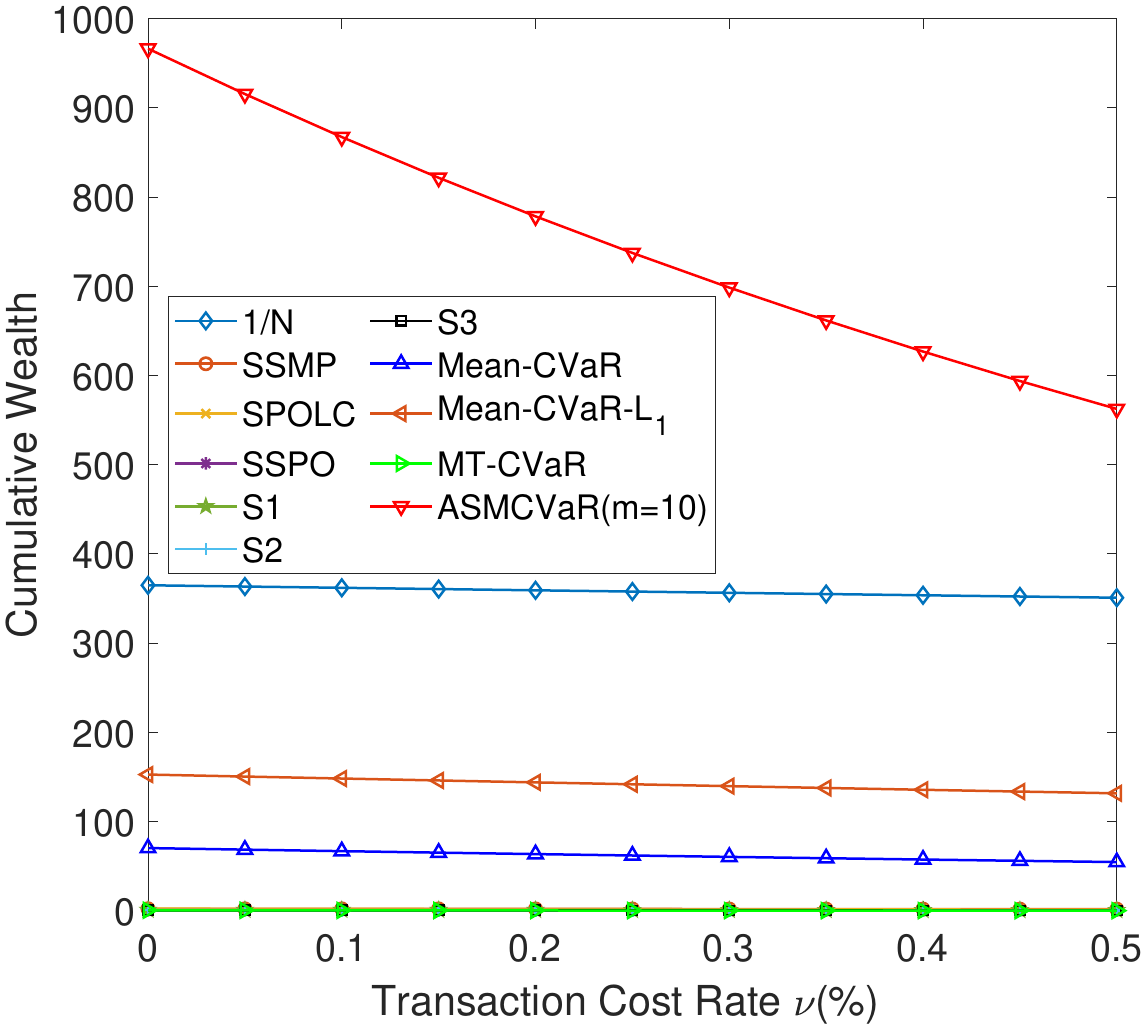}}
\subfloat[FF100MEOP]{
\centering
\includegraphics[width=0.42\textwidth]{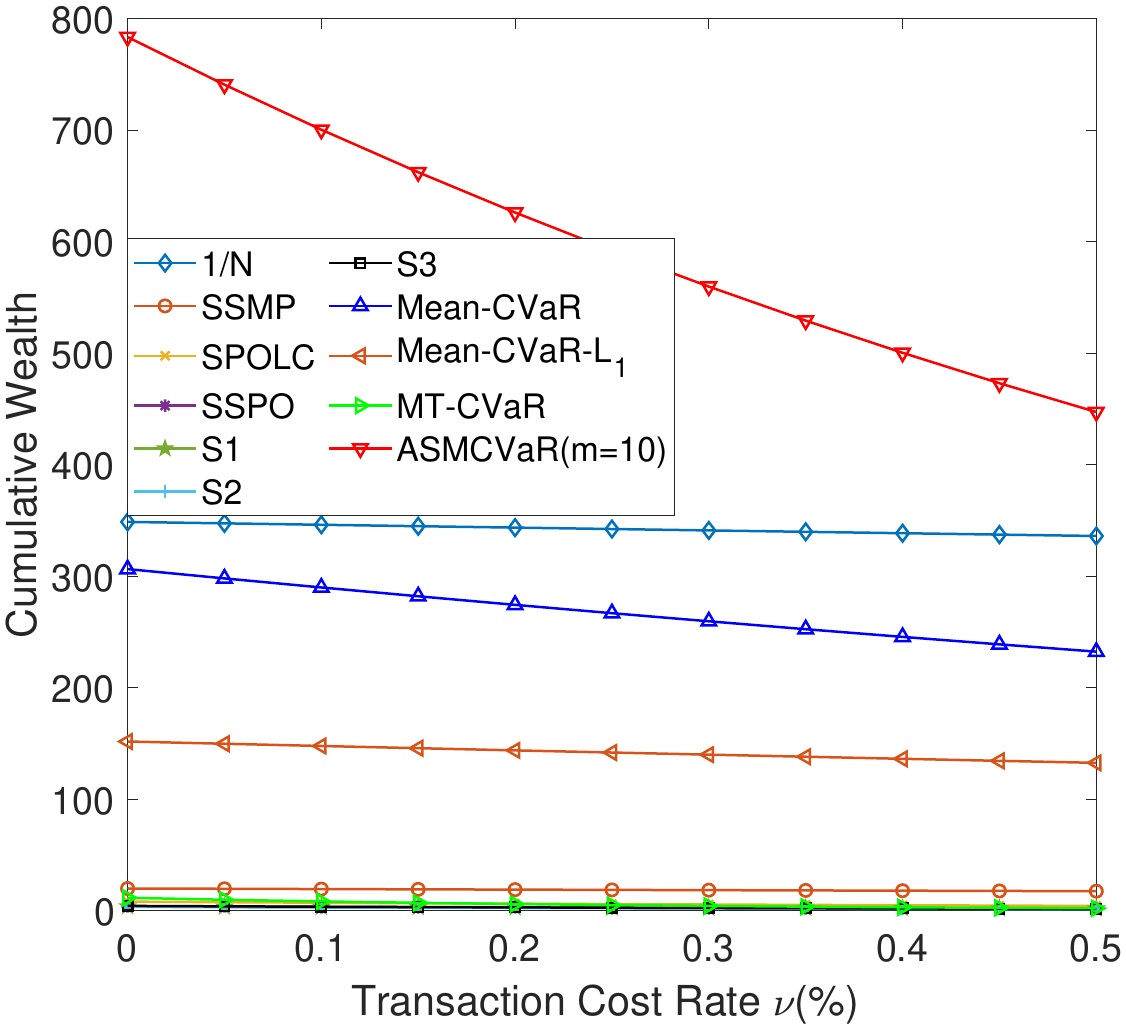}}
\caption{Final cumulative wealths of different portfolio optimization models w.r.t. transaction cost rate $\nu$ on $6$ benchmark data sets.}
\label{fig:tc}
\end{figure*}

\subsection{Transaction Cost Plots (Figure \ref{fig:tc})}
\label{app:tc}

\end{document}